\documentclass[aop,MSNbibl,dvips]{arximspdf}

\usepackage{mathbh}

\doi{10.1214/13-AOP882} 
\volume{42}
\issue{4}
\pubyear{2014}
\firstpage{1666}
\lastpage{1698}

\makeatletter
\newcommand{\varPsi}{\Psi}
\newcommand{\varTheta}{\Theta}
\newcommand{\varPi}{\Pi}
\newcommand{\varSigma}{\Sigma}
\newcommand{\varGamma}{\Gamma}
\newcommand{\eqref}[1]{(\ref{#1})}
\renewcommand{\backslash}{\setminus}
\renewcommand{\scriptscriptstyle}{}
\renewcommand{\phi}{\varphi}
\newcommand{\id}{\one}
\newcommand{\ssup}[1]{{\scriptscriptstyle{({#1}})}}
\newcommand{\one}{\mathbh{1}}
\newcommand{\Prob}{\operatorname{Prob}}
\newcommand{\Deltad}{\Delta}
\newcommand{\dd}{d}
\newcommand{\R}{\mathbb{R}}
\newcommand{\Z}{\mathbb{Z}}
\newcommand{\N}{\mathbb{N}}
\newcommand{\E}{\mathbb{E}}
\renewcommand{\P}{\mathbb{P}}
\newcommand{\1}{\mathbh{1}}

\newtheorem{theorem}{Theorem}[section]
\newproclaim{remark}{Remark}
\newtheorem{lemma}[theorem]{Lemma}
\newtheorem{prop}[theorem]{Proposition}
\makeatother

\begin{document}
\begin{frontmatter}

\title{Localisation and ageing in the parabolic Anderson model with
Weibull potential}
\runtitle{Localisation and ageing in the PAM with Weibull potential}

\begin{aug}
\author[A]{\fnms{Nadia} \snm{Sidorova}\corref{}\ead[label=e1]{n.sidorova@ucl.ac.uk}}
\and
\author[A]{\fnms{Aleksander} \snm{Twarowski}\ead[label=e2]{a.twarowski@ucl.ac.uk}}
\runauthor{N. Sidorova and A. Twarowski}
\affiliation{University College London}
\address[A]{Department of Mathematics\\
University College London\\
Gower Street\\
London WC1 E6BT\\
United Kingdom\\
\printead{e1}
\\
\phantom{E-mail: }\printead*{e2}} 
\end{aug}

\received{\smonth{6} \syear{2012}}
\revised{\smonth{4} \syear{2013}}

%
\begin{abstract}
The parabolic Anderson model is the Cauchy problem for the heat
equation on the integer lattice with a random potential $\xi$.
We consider the case when $\{\xi(z)\dvtx z\in\mathbb{Z}^d\}$
is a collection of independent identically distributed random variables
with Weibull distribution with parameter $0<\gamma<2$,
and we assume that the solution is initially localised in the origin.
We prove that, as time goes to infinity, the solution completely localises
at just one point with high probability, and we identify the asymptotic
behaviour of the localisation site.
We also show that the intervals between the times when the solution
relocalises from one site to another
increase linearly over time, a phenomenon known as ageing.
\end{abstract}

%

\end{frontmatter}

\section{Introduction and main results}\label{intro}

\subsection{Parabolic Anderson model}

We consider the heat equation with random potential on the integer
lattice $\Z^d$ and study
the Cauchy problem with localised initial condition,
%
%
\begin{equation}
\label{pam} %
\begin{array} {@{}rcl@{}}
\partial_t u(t,z) & = & \Deltad u(t,z)+\xi(z)u(t,z), \qquad(t,z)\in(0,
\infty)\times \Z^d,
\\[5pt]
u(0,z) & = & \one_{\{0\}}(z), \qquad z\in\Z^d, \end{array}
\end{equation}
where
\[
(\Deltad f) (z)=\sum_{y\sim z} \bigl[f(y)-f(z) \bigr],
\qquad z\in\Z^d, f\dvtx\Z^d\to\R
\]
is the discrete Laplacian, and the potential
$ \{\xi(z) \dvtx z\in\Z^d \}$ is a collection
of independent identically distributed random variables.
The problem \eqref{pam} and its variants are often called the
\emph{parabolic Anderson model}.

The model originates from the seminal work \cite{An58} of the Nobel
laureate P.~W. Anderson, who used the Hamiltonian $\Delta+\xi$ to
describe electron localisation inside a semiconductor, a phenomenon now
known as Anderson localisation. The parabolic version of the model
appears naturally in the context of reaction--diffusion equations; see
\cite{CM94,M94},
describing a system of noninteracting particles diffusing
in space according to the Laplacian $\Delta$ and branching at rate
$\xi
(z)\,dt$ at any given point~$z$. It turns out that
the solution $u(t,z)$ gives the average number of such particles at
time $t$ at location $z$.

\subsection{Intermittency and localisation}

A lot of mathematical attention to the parabolic Anderson model over
the last 30 years has been due to the fact that it exhibits
the \emph{intermittency
effect}. In general, a random model is said to be intermittent if its
long-term behaviour cannot be described
using an averaging principle; see \cite{ZM87}.
In the context of the parabolic Anderson model,
this means that, for large times $t$, the solution $u(t,z)$ is mainly
concentrated on a small number of remote random islands;
see \cite{GK05} for a survey.

The long-term behaviour of the parabolic Anderson model is determined
by the upper tail of
the underlying distribution of the potential $\xi$, and it is believed
that the intermittency is more pronounced
for heavier tails. However, an initial approach to understanding
intermittency was proposed for light-tailed
potentials (those with finite exponential moments). It was suggested to study
large time asymptotics of the moments of the total mass of the solution
\[
U(t)=\sum_{z\in\Z^d} u(t,z),
\]
which are finite for such potentials.
The model was defined as intermittent if higher moments exhibited a
faster growth rate, and it was proved in \cite{GM90}
that the parabolic Anderson model is intermittent in this sense.
This method, however, does not work for heavy-tailed potentials
(those with infinite exponential moments), as for them the moments of
$U(t)$ are infinite.
Such distributions include the exponential distribution and all
heavier-tailed distributions.

In order to understand the intermittent picture in more detail, it
proved to be useful to
study various large-time asymptotics of the total mass $U(t)$, as they
provided some insight into the geometry of
the intermittent islands. It was shown in \cite{HKM06} that there are
four types of behaviour the parabolic Anderson
model can exhibit depending on the tail of the underlying distribution.
The prime examples from each class are the
following distributions:
\begin{enumerate}[(3)]
\item[(1)] Weibull distribution with parameter $\gamma>1$, that is,
$F(x)=1-e^{-x^{\gamma}}$.
\item[(2)] Double-exponential distribution with parameter $\rho>0$,
that is, $F(x)=1-e^{-e^{x/\rho}}$.
\item[(3)] ``Almost bounded'' distributions, including some unbounded
distributions with tails lighter than double-exponential
and some bounded distributions.
\item[(4)] Other bounded distributions.
\end{enumerate}
The asymptotics of the total mass $U(t)$ was studied in \cite{GM98} for
cases (1) and (2), in \cite{HKM06}
for case (3) and in \cite{BK} for case (4).
Heuristics based on the asymptotics of $U(t)$ suggests that the
intermittent islands will be single lattice points
in case (1), bounded regions in case (2) and of size growing to
infinity in cases (3)
and~(4). However, a rigorous geometric picture of intermittency has not
been well understood.
In particular, it is not clear how many intermittent islands are needed
to carry the total mass of the solution,
and where those islands are located.

Moreover, the four classes above only cover light-tailed potential, and
the class
of all heavy-tailed distributions should be included to
complete the picture. The prime examples of such distributions are
\begin{enumerate}[(0b)]
\item[(0a)] Pareto distributions, that is, $F(x)=1-x^{-\alpha}$,
$\alpha>d$;
\item[(0b)] Weibull potentials with parameter $\gamma\le1$.
\end{enumerate}
Heavy-tailed potentials were first studied in \cite{HMS08}, and it
turned out that
the asymptotics of $U(t)$ in this case becomes nondeterministic and
difficult to control.
It was suggested to study the nondeterministic nature
of $U(t)$ using extreme value theory and point processes techniques.
This approach was further developed
in \cite{KLMS09}, where the intermittency was fully described in its
original geometric sense for Pareto potentials (0a).
Polynomial tails are the heaviest tails for which the solution of the
parabolic Anderson model still exists (see \cite{GM90}),
and one expected the localisation islands to be small and not numerous.
It was proved that the extreme form of this conjecture
is true, namely, that there is only one localisation island consisting
of only one site. In other words, at any time
the solution is localised at just one point with high probability, a
phenomenon called \emph{complete localisation}.

It is a challenging problem to describe geometric intermittency for
lighter tails.
In \cite{GKM06}, intermittent islands were described for potentials
from classes (1) and (2), but
the question about the number of islands remained open.
Case (0b) was studied in \cite{LM},
and it was shown that the solution is localised on an island of size
$o(\frac{t (\log t)^{1/\gamma-1}}{\log\log t})$. However,
it was believed that a much smaller region should actually contribute
to the solution.

In this paper, we
assume that the potential has \emph{Weibull distribution} with
parameter $\gamma>0$, that is, the distribution function of each $\xi
(z)$ is
%
%
\begin{equation}
\label{wei} F(x)=\Prob \bigl\{\xi(z)<x \bigr\}=1-e^{-x^{\gamma
}},\qquad x
\ge0.
\end{equation}
We focus on $0<\gamma<2$, which covers case (0b) and partly case (1).
We prove
that for such potentials the solution of the parabolic Anderson model
completely localises
at just one single site, exhibiting the strongest form of intermittency
similar to the Pareto case (0a).
This was plausible for $0<\gamma<1$ as in this case the spectral gap of
the Anderson Hamiltonian
$\Delta+\xi$ in a relevant $t$-dependent large box tends to infinity,
but is quite surprising for
the exponential distribution ($\gamma=1$) where the spectral gap is
bounded, and even more so
for $1<\gamma<2$ where the spectral gap tends to zero.
We identify the localisation site explicitly in terms of the potential
$\xi$ and describe its scaling limit.


For all sufficiently large $t$ (so that $\log\log t$ is well defined), denote
%
%
\begin{equation}
\label{phi} \varPsi_t(z)=\xi(z)-\frac{|z|}{\gamma t}\log\log t,\qquad
z\in\Z^d,
\end{equation}
and let $Z_t^{\ssup1}$ be such that
\[
\varPsi_t\bigl(Z_t^{\ssup1}\bigr)=\max
_{z\in\Z^d} \varPsi_t(z).
\]
The existence of $Z_t^{\ssup1}$ will be proved in Lemma~\ref{welldef}.

Denote by $|x|$ the $\ell^1$-norm of $x\in\R^d$, and denote by
$\Longrightarrow$ weak convergence.

%
\begin{theorem}[(Complete localisation)]\label{main_w}
Let $0<\gamma<2$. As $t\to\infty$,
\[
\lim_{t\to\infty}\frac{u(t,Z_t^{\ssup1})}{U(t)}= 1 \qquad\mbox {in probability}.
\]
\end{theorem}
%

%
\begin{remark}
It is easy to see that the solution cannot be localised
at one point for all large times
$t$ since occasionally it has to relocalise continuously from one site
to another, and at those periods the solution
will be concentrated at more than one point. It was shown in \cite
{KLMS09} that for Pareto
potentials the solution in fact remains localised at just two points at
all large times $t$ almost surely. We conjecture that
the same is true for Weibull potentials with $0<\gamma<2$.
\end{remark}
%

%
\begin{remark}
There is a chance that our proof could be adjusted to
the case \mbox{$\gamma=2$}. However,
new ideas are required to deal with $\gamma>2$, and there is a high
chance that complete localisation will
simply fail in that case. The technical reasons why our proof breaks
down for $\gamma\ge2$ are explained
in Remark~\ref{rem1} and Remark~\ref{rem2} in Sections~\ref{s_neg}
and~\ref{s_loc}, respectively.
\end{remark}
%

%
\begin{theorem}[(Scaling limit for the localisation site)]\label{main_z}
Let $\gamma>0$. Then
%
\[
\frac{Z_t^{\ssup1}}{r_t} \Longrightarrow X^{\ssup1},
\]
as $t\to\infty$ where
%
%
\begin{equation}
\label{rrr} r_t=\frac{t(\log t)^{1/\gamma-1}}{\log\log t}\vadjust{\goodbreak}
\end{equation}
and
$X^{\ssup1}$ is an $\R^d$-valued random variable with independent
exponentially distributed coordinates
with parameter $d^{1-1/\gamma}$ and uniform random signs, that is,
with density
\[
p^{\ssup1}(x) =\frac{d^{d(1-1/\gamma)}}{2^d}\exp \bigl\{ -d^{1-1/\gamma
}|x| \bigr\},
\qquad x\in\R^d.
\]
\end{theorem}
%

%
\begin{remark}
Although we prove Theorem~\ref{main_z} for all $\gamma
>0$, it only describes the scaling
limit for the concentration site for $0<\gamma<2$ as otherwise the
solution may not be localised
at $Z_t^{\ssup1}$.
\end{remark}
%

%
\begin{remark}
This scaling limit agrees with the scaling limit for the
centre of the intermittent island obtained in \cite{LM}
for $0<\gamma\le1$.
However, according to Theorem~\ref{main_w}, this island is now of
radius zero (being a single point) rather than $o(r_t)$,
and the result holds for the wider range $0<\gamma<2$.
\end{remark}
%

\subsection{Ageing}

The notion of \emph{ageing} is a key paradigm in studying the long-term
dynamics of large disordered systems.
A system exhibits
ageing if, being in a certain state at time $t$, it is likely to remain
in this state for some time $s(t)$ which
depends increasingly, and often linearly, on the time $t$. Roughly
speaking, the system becomes increasingly more
conservative and reluctant to change.

The ageing phenomenon has been extensively studied for disordered
systems such as trap models and spin glasses;
see \cite{BC} and references therein.
In the context of the parabolic Anderson model, a certain form of
ageing based on correlations was studied for some
time-dependent potentials in \cite{AD,DD}, and it was shown that such
systems exhibit no ageing.
The recent paper \cite{GS} dealt with potentials from class (1) and
studied the
correlation ageing (which gives only indirect information about the
evolution of localisation) and more explicit
annealed ageing (which, in contrast to the quenched setting, is based
on the evolution
of the islands contributing to the solution averaged over the environment).
It was shown that these two forms of ageing are similar, and
somewhat surprisingly, ageing was observed for Weibull potentials with
parameter $\gamma>2$ but not for
heavier-tailed Weibull potentials with parameter $1< \gamma\le2$.

The explicit ageing in the quenched setting has so far only been
observed for Pareto potentials; see \cite{MOS}. In that case, the
solution completely localises at just one point
and ageing of the parabolic Anderson model is equivalent to ageing of
the concentration site process.
In this paper, we use a similar approach to show that
the parabolic Anderson model with Weibull potential with parameter
$0<\gamma<2$ exhibits ageing
as well. Notice that, remarkably, this is in sharp contrast to the
absence of annealed and correlation ageing
observed for $\gamma>1$ in \cite{GS}.

For each $t>0$, denote
\[
T_t=\inf \bigl\{s>0\dvtx Z_{t+s}^{\ssup1}\neq
Z_t^{\ssup1} \bigr\}.
\]
%
%
\begin{theorem}[(Ageing)]
\label{main_a}
Let $\gamma>0$.
As $t\to\infty$
\[
\frac{T_t}{t}\Longrightarrow\varTheta,
\]
where $\varTheta$ is a nondegenerate almost surely positive random variable.
\end{theorem}
%

%
\begin{remark}
In the proof of Theorem~\ref{main_a}, we identify the
distribution function of $\varTheta$
as a certain integral over $\R^d\times\R$.
\end{remark}
%

%
\begin{remark}
Although we prove Theorem~\ref{main_a} for all $\gamma
>0$, it only characterises the ageing
behaviour of the parabolic Anderson model for $0<\gamma<2$ as otherwise
the solution may not be localised
at $Z_t^{\ssup1}$.
\end{remark}

\subsection{Outline of the proofs}

It follows from \cite{GM90}, Theorem~2.1, that
the parabolic Anderson model with Weibull potential
possesses a unique nonnegative solution $u\dvtx(0,\infty)\times\Z
^d\to[0,\infty)$,
which has a \emph{Feynman--Kac representation}
\[
u(t,z)=\E_0 \biggl[\exp \biggl\{\int_0^t
\xi(X_s)\,\dd s \biggr\}\id\{ {X_t=z\}} \biggr],\qquad
(t, z)\in(0,\infty)\times\Z^d,
\]
where $(X_s \dvtx s\ge0)$ is a
continuous-time simple random walk on the lattice $\Z^d$ with generator
$\Deltad$, and $\P_z$ and $\E_z$ denote
the corresponding probability and expectation given that the random
walk starts at $z\in\Z^d$.

The Feynman--Kac formula suggests that the main contribution to the
solution $u$ at time $t$ comes from paths $(X_s)$
spending a lot of time at sites $z$ where the value $\xi(z)$ of the
potential is high but which are
reasonably close to the origin so that the random walk would have a
fair chance of reaching them in time $t$.
It turns out that the functional $\varPsi_t$ defined in \eqref{phi}
captures this trade-off, being the difference of the energetic term
$\xi
(z)$ and an entropic term
responsible for the cost of going to a point $z$ in time $t$ and
staying there.
Furthermore, the maximiser $Z_t^{\ssup1}$ of $\varPsi_t$ turns out to
be the site
where the solution $u$ is localised at time $t$.

In order to prove this, we decompose the solution $u$ into the sum
\[
{u(t,z)=u_1(t,z)+u_2(t,z)}
\]
according to two groups of paths ending at $z$:
\begin{enumerate}[(II)]
\item[(I)] paths visiting $Z_t^{\ssup1}$ before time $t$ and staying
in the ball
$B_t$ centred in the origin with radius
$|Z_t^{\ssup1}|(1+\rho_t)$, where $\rho_t$ is a certain function
tending to zero;
\item[(II)] all other paths.
\end{enumerate}
We show that $u_1$ localises around $Z_t^{\ssup1}$ and that the total
mass of
$u_2$ is negligible.

To prove the localisation of $u_1$, we use spectral analysis of the
Anderson Hamiltonian $\Delta+\xi$
in the ball $B_t$. In order to do so, we show that, although the
spectral gap tends to zero for $\gamma>1$,
it is still reasonably large. We suggest a new technique which allows
us to show that
the principal eigenfunction just manages to localise at $Z_t^{\ssup1}$.
Then we use a result from \cite{GKM06} to show that this is sufficient
for the localisation of $u_1$.

In order to prove that the total mass of $u_2$ is negligible, we notice
that the paths from the second group
fall into one of the following three subgroups:
\begin{enumerate}[(3)]
\item[(1)] paths having the maximum of the potential at the point
$Z_t^{\ssup1}$
but making more than $|Z_t^{\ssup1}|(1+\rho_t)$ steps;
\item[(2)] paths having the maximum of the potential not at the point
$Z_t^{\ssup1}$, with the maximum being reasonably large;
\item[(3)] paths missing all high values of the potential.
\end{enumerate}

In Section~\ref{s_neg}, we show that the total mass of the paths
corresponding to each group is negligible.
In all cases, this is due to an imbalance between
the energetic forces (which do not contribute enough if the site
$Z_t^{\ssup1}$
is not visited) and entropic forces
(as the probabilistic cost is too high if a path is too long), as well
as to the fact that the gap between
$\varPsi_t(Z_t^{\ssup1})$ and the second largest value of $\varPsi
_t$ is too large.

Denote by
$Z_t^{\ssup2}$ a point where the second largest value of $\varPsi_t$
is attained,
that is,
\[
\varPsi_t \bigl(Z_t^{\ssup2} \bigr)=
\max \bigl\{\varPsi_t(z)\dvtx z\in\Z^d, z\neq
Z_t^{\ssup1} \bigr\}.
\]
In order to find the scale of growth of $\varPsi_t(Z_t^{\ssup
1})-\varPsi_t(Z_t^{\ssup2})$ as
well as of $Z_t^{\ssup1}$ and $Z_t^{\ssup2}$ we extend the
point processes techniques developed in \cite{HMS08} and \cite{KLMS09}.
For sufficiently large $t$, we denote
\[
a_t=(d\log t)^{1/\gamma} \quad\mbox{and}\quad d_t=(d
\log t)^{1/\gamma-1}.
\]
Further, for all $z\in\Z^d$ and all sufficiently large $t$, we denote
%
%
\begin{equation}
\label{yy} Y_{t,z}=\frac{\varPsi_t(z)-a_{r_t}}{d_{r_t}},
\end{equation}
where $r_t$ is defined by \eqref{rrr},
and define a point process
%
%
\begin{equation}
\label{pp} \varPi_t=\sum_{z\in\Z^d}
\varepsilon_{(zr_t^{-1},Y_{t,z})},
\end{equation}
where we write $\varepsilon_x$ for the Dirac measure in $x$. In
Section~\ref{s_ppp}, we show that
the point processes $\varPi_t$ are well defined on a carefully chosen
domain, and that they converge in law to a Poisson point process with
certain density.
This allows us to analyse the joint distribution of the random
variables $Z_t^{\ssup1}$, $Z_t^{\ssup2}$, $\varPsi_t(Z_t^{\ssup
1})$, $\varPsi_t(Z_t^{\ssup2})$
and, in particular, prove Theorem~\ref{main_z}.

Finally, to prove ageing, we argue that due to the form of the
functional $\varPsi_t$
the probability of $\{Z^{\ssup1}_{t+wt}= Z_t^{\ssup1}\}$, for each
$w>0$, is
roughly equal to
%
%
\begin{equation}
\label{inti} \int_{\R^d\times\R}\Prob \bigl\{\varPi_t(
\dd x\times\dd y)=1, \varPi_t \bigl(D_{w}(x,y) \bigr)=0
\bigr\},
\end{equation}
where
%
%
\begin{eqnarray}
\label{dc} D_{w}(x,y)&=& \biggl\{(\bar x,\bar y)\in
\R^d\times\R\dvtx y+\frac{w\theta
|x|}{1+w}\le\bar y+\frac{w\theta|\bar x|}{1+w}
\biggr\}
\nonumber
\\[-8pt]
\\[-8pt]
&&{}\cup \bigl(\R^d\times[y,\infty ) \bigr),
\nonumber
\end{eqnarray}
and
%
%
\begin{equation}
\label{thth} \theta=\gamma^{-1}d^{1-1/\gamma}.
\end{equation}
In particular, the integral in \eqref{inti} converges to the
corresponding finite integral with respect to
the Poisson point process $\varPi$ as $t\to\infty$. This proves Theorem~\ref{main_a} since that integral is a continuous function
of $w$ decreasing from
one to zero as $w$ varies from zero to infinity and so it is the tail
of a distribution function.

The paper is organised as follows. In Section~\ref{s_pre}, we introduce
notation and prove some
preliminary results. In Section~\ref{s_ppp}, we develop a point
processes approach, analyse
the joint distribution of $Z_t^{\ssup1}$, $Z_t^{\ssup2}$, $\varPsi
_t(Z_t^{\ssup1})$, $\varPsi_t(Z_t^{\ssup2})$
and prove Theorem~\ref{main_z}.
In Section~\ref{s_neg}, we deal with the total mass corresponding to
the paths from groups (1)--(3) and show
that it is negligible. In Section~\ref{s_loc}, we discuss the
localisation of $u_1$ and prove Theorem~\ref{main_w}.
Finally, in Section~\ref{s_age}, we study ageing and prove Theorem~\ref{main_a}.


\section{Preliminaries}
\label{s_pre}

We focus on potentials with Weibull distribution \eqref{wei} with
parameter $0<\gamma<2$.
However, most of our point processes results can be obtained for all
$\gamma>0$ at no additional cost. Therefore,
we will assume $\gamma>0$ in Sections~\ref{s_pre}, \ref{s_ppp}
and~\ref{s_age},
and restrict ourselves to the case
$0<\gamma<2$ in Sections~\ref{s_neg} and~\ref{s_loc}.

\subsection{Extreme value notation and preliminary results}

We denote the upper order statistics of the potential $\xi$ in the
centred ball
of radius $r>0$ by
\[
\xi_r^{\ssup1}=\max_{|z|\le r}\xi(z)
\]
and
\[
\xi^{\ssup i}_r=\max \bigl\{\xi(z)\dvtx|z|\le r, \xi(z)< \xi
_r^{\ssup
{i-1}} \bigr\}
\]
for $2\le i\le\ell_r$, where $\ell_r$ is the number of points in the
ball. Observe that throughout the paper
we use the $\ell^1$-norm.

Let $0<\rho<\sigma<1/2$ and for all sufficiently large $r$ let
\begin{eqnarray*}
F_r&=& \bigl\{z\in\Z^d\dvtx|z|\le r, \exists i\le
r^{\rho} \mbox{ such that }\xi(z)=\xi_{r}^{\ssup i}
\bigr\},
\\
G_r&=& \bigl\{z\in\Z^d\dvtx|z|\le r, \exists i\le
r^{\sigma} \mbox{ such that }\xi(z)=\xi_{r}^{\ssup i}
\bigr\}.
\end{eqnarray*}
The sets $F_r$ and $G_r$ contain the sites in the centred ball of
radius $r$
where the highest $\lfloor r^{\rho}\rfloor$ and $\lfloor r^{\sigma
}\rfloor$ values of the potential $\xi$
are achieved, respectively.


%
\begin{lemma}
\label{l_asasas}
Almost surely
\[
\xi_{r}^{\ssup1}\sim(d\log r)^{1/\gamma} \qquad\mbox{as }
r\to\infty.
\]
\end{lemma}


\begin{pf}
This result was proved in \cite{HMS08} for the case
$0<\gamma\le1$ but it can be
easily extended to all $\gamma>0$ by observing that $\zeta(z)=\xi
(z)^{\gamma}, z\in\Z$,
are exponential identically distributed random variables.
Denote the maximum of the potential $\zeta$ by
\[
\zeta_r^{\ssup1}=\max_{|z|\le r}\zeta(z).
\]
%
Since $\xi_r^{\ssup1}= (\zeta_r^{\ssup1} )^{1/\gamma}$ and
$\zeta_r^{\ssup1}\sim d\log r$
by \cite{HMS08}, Lemma~4.1, with $\gamma=1$, we obtain
the required asymptotics.
\end{pf}

For all $c\in\R$, $z\in\Z^d$, and all sufficiently large $t$ define
\[
\varPsi_{t,c}(z)=\varPsi_t(z)+\frac{c|z|}{t}.
\]
Denote by $Z_t^{\ssup{1,c}}$ and $Z_t^{\ssup{2,c}}$ points where the
first and second largest values of the
functional $\varPsi_{t,c}$ are achieved, that is,
%
%
\begin{equation}
\label{maxc} %
\begin{array} {@{}rcl@{}}
\varPsi_{t,c} \bigl(Z_t^{\ssup{1,c}} \bigr)&=&\max \bigl\{
\varPsi_{t,c}(z) \dvtx z\in\Z^d \bigr\},
\\[5pt]
\varPsi_{t,c} \bigl(Z_t^{\ssup{2,c}} \bigr)&=&\max \bigl
\{ \varPsi_{t,c}(z)\dvtx z\in\Z^d, z\neq
Z_t^{\ssup{1,c}} \bigr\}. \end{array} %
\end{equation}
Observe that $\varPsi_{t}=\varPsi_{t,0}$ and so $Z_t^{\ssup
1}=Z_t^{\ssup{1,0}}$ and
$Z_t^{\ssup2}=Z_t^{\ssup{2,0}}$.
We are mostly interested in the case $c=0$, but some
understanding of the general case is needed for Lemma~\ref{l_u2}. This
is explained more carefully in Remark~\ref{rem3}
in Section~\ref{s_ppp}.

%
\begin{lemma}
\label{welldef}
For each $c$, the maximisers $Z_t^{\ssup{1,c}}$ and $Z_t^{\ssup{2,c}}$
(and, in particular, $Z_t^{\ssup1}$ and $Z_t^{\ssup2}$) are well
defined for all
sufficiently large $t$
almost surely.
\end{lemma}

\begin{pf} Observe that $\varPsi_{t,c}(0)>0$ and $\varPsi_{t,c}(1)>0$
almost surely if $t$ is large enough.
On the other hand, by Lemma~\ref{l_asasas} for all sufficiently large
$t$ there
exists a random radius $\rho(t)>0$ such that, almost surely,
\[
\xi(z)\le\xi_{|z|}^{\ssup1}\le\bigl(2d
\log|z|\bigr)^{1/\gamma} \le\frac{|z|}{\gamma t}\log\log t-\frac
{c|z|}{t}\qquad
\mbox{for all }|z|>\rho(t).
\]
Hence, $\varPsi_{t,c}(z)\le0$ for all $|z|>\rho(t)$ and so
$\varPsi_{t,c}$ takes only finitely many positive values.
This implies that the maxima
in \eqref{maxc} exist for all $c$. The existence of $Z_t^{\ssup1}$
and $Z_t^{\ssup2}$
follows as a particular case when $c=0$.
\end{pf}

Choose
\[
\cases{ %
\beta\in(1-1/\gamma,1/\gamma)&\quad$\mbox{if }1
\le\gamma<2$,
\cr
\beta=0 &\quad$\mbox{if }0< \gamma<1$. 
}
\]
Observe that $\beta\ge0$
and define
%
%
\begin{equation}
\label{mu} \mu_r=(\log r)^{-\beta}
\end{equation}
for all $r$ large enough.
For $0< \gamma<1$, the gaps between higher order statistics of the
potential get larger (as $r\to\infty$) and
the auxiliary scaling function $\mu_r$ is not needed (so that we can
simply set $\mu_r=1$ as above). For $\gamma=1$, the gaps are of finite
order, and for $\gamma>1$ they tend to zero, and an extra
effort is required to control this effect. This is done by the
correction term $\mu_r$.
It is essential for the choice of $\mu_r$ that, on the one hand, it is
negligible with respect to $d_r$
and so with respect to the gap $\varPsi_t(Z_t^{\ssup1})-\varPsi
_t(Z_t^{\ssup2})$ (which is
achieved by the condition $\beta>1-1/\gamma$)
and on the other hand $-\log\mu_r$ must be smaller than $\log\xi
_r^{\ssup1}$ (which is guaranteed by $\beta<1/\gamma$).
However, this method only works
for $\gamma<2$ as the interval $(-1/\gamma+1,1/\gamma)$ is empty
otherwise. This is explained in more detail
in Remark~\ref{rem1} in Section~\ref{s_neg}.

We introduce four auxiliary positive scaling functions
$f_t\to0$, $g_t\to\infty$, $\lambda_t\to0$, $\rho_t\to0$ satisfying
the following conditions as $t\to\infty$:
%
%
\begin{eqnarray}
\label{a} &&\mbox{(a)} \quad f_t^{-1},
g_t, \lambda_t^{-1}, \rho_t^{-1}
\mbox{ are } o(\log\log t),
\\
\label{b} &&\mbox{(b)} \quad g_t\rho_t
\lambda_t^{-1}\to0. 
\end{eqnarray}
%


Further, we define
\[
k_t= \bigl\lfloor(r_tg_t)^{\rho}
\bigr\rfloor\quad\mbox{and}\quad m_t= \bigl\lfloor(r_tg_t)^{\sigma}
\bigr\rfloor.
\]


For any $c\in\R$, we introduce the event
%
%
\begin{eqnarray}
\label{eee} &&\mathcal{E}_c(t)= \bigl\{r_tf_t<\bigl |
Z_t^{\ssup1}\bigr |<r_tg_t,
\varPsi_t\bigl(Z_t^{\ssup1}\bigr)-\varPsi
_t\bigl(Z_t^{\ssup2}\bigr)>d_t
\lambda_t,
\nonumber
\\
&&\phantom{\mathcal{E}_c(t)= \bigl\{} \varPsi_t\bigl(
Z_t^{\ssup1}\bigr)>a_{r_t}-d_t
g_t, \varPsi_t\bigl( Z_t^{\ssup2}
\bigr)>a_{r_t}-d_t g_t,
\\
&&\phantom{\hspace*{128pt}}\bigl |Z^{\ssup{1,c}}_t\bigr |<r_tg_t,
\bigl |Z^{\ssup{2,c}}_t\bigr |<r_tg_t \bigr\}.
\nonumber
\end{eqnarray}

For any $x,y\in\R$, we denote by $x \wedge y$ and $x\vee y$ the
minimum and the maximum of $x$ and $y$,
respectively, and we denote $x_{-}=-x\vee0$.


\subsection{Geometric paths on the lattice}
\label{s_geo}

For each $n\in\N\cup\{0\}$ denote by
\[
\mathcal{P}_n= \bigl\{y=(y_0,\ldots,y_n)\in
\bigl(\Z^d \bigr)^{n+1}\dvtx|y_i-y_{i-1}|=1
\mbox{ for all }1\le i\le n \bigr\}
\]
the set of all geometric paths in $\Z^d$.
Define
\[
q(y)=\max_{0\le i\le n}\xi(y_i) \quad\mbox{and}\quad
p(y)=\max_{0\le i\le n} |y_i-y_0|,
\]
and denote by $z(y)$ a point $y_i$ of the path $y$ such that $\xi(y_i)=q(y)$.

Let $(\tau_i)$, $i\ge0$, be waiting times
of the random walk $(X_s)$, which are independent exponentially
distributed random variables with
parameter $2d$. Denote by $\mathsf{E}$ the expectation with respect to
$(\tau_i)$.
For each $y\in\mathcal{P}_n$, denote by
\begin{eqnarray*}
&&P(t,y)= \{X_0=y_0, X_{\tau_0+\cdots+\tau_{i-1}}=y_i
\mbox{ for all }1\le i\le n,
\\
&&\phantom{\hspace*{98pt}} \mbox{and } t-\tau_n\le
\tau_0+\cdots+\tau_{n-1}< t \}
\end{eqnarray*}
the event that the random walk has the trajectory $y$ up to time $t$. Here,
we assume that the random walk is continuous from the right. Denote by
%
%
\begin{equation}
\label{uty} U(t,y)=\E_0 \biggl[\exp \biggl\{\int
_0^t\xi(X_s)\,\dd s \biggr\}\1
_{P(t,y)} \biggr]
\end{equation}
the contribution of the event $P(t,y)$ to the total mass of the
solution $u$ of the parabolic Anderson model.
%
%

For any set $A\subset\Z^d$ and any geometric path $y\in\mathcal
{P}_n$ denote
\[
n_{+}(y,A)=\bigl | \{0\le i\le n\dvtx y_i\in A \}\bigr | \quad
\mbox{and}\quad n_{-}(y,A)=\bigl | \{0\le i\le n\dvtx y_i\notin
A \}\bigr |.
\]
We call a set $A\subset\Z^d$ totally disconnected if $|x-y|\neq1$
whenever $x,y\in A$.

%
\begin{lemma}
\label{l_count}
Let $A$ be a totally disconnected finite subset of $\Z^d$,
and $y\in\mathcal{P}_n$ for some $n$.
Then
\[
n_+(y,A)\le\frac{n-p(y)}{2}+|A|\wedge \biggl\lceil\frac
{p(y)+1}{2} \biggr
\rceil. 
\]
\end{lemma}

\begin{pf}
Let $i(y)=\min\{i\dvtx|y_i-y_0|=p(y)\}$ and denote $z=y_{i(y)}$. Similarly
to \cite{KLMS09}, page 371, we first
erase loops that the path $y$ may have made
before reaching $z$ for the first time and extract from $(y_0,\ldots
,y_{i(y)})$ a self-avoiding
path $(y_{i_0},\ldots,y_{i_{p(y)}})$ starting at $y_0$ of length $p(y)$,
where we take $i_0=0$ and
\[
i_{j+1}=\min \bigl\{i\dvtx y_l\neq y_{i_j}\
\forall l\in \bigl[i,i(y) \bigr] \bigr\}.
\]
%
Since this path is self-avoiding and has length $p(y)$,
at most $|A|\wedge\lceil\frac{p(y)+1}{2}\rceil$ of its points belong
to $A$.
Next, for each $0\le j\le p(y)-1$, we consider the path
$(y_{i_j+1},\ldots,y_{i_{j+1}-1})$, which was
removed during erasing the $j$th loop. It contains an even number
$i_{j+1}-{i_j}-1$
of steps
and at most half of them
belong to $A$ since $A$ is totally disconnected.
Finally, the remaining piece $(y_{i_{p(y)}+1},\ldots,y_n)$ consists of
$n-i_{p(y)}$ points, and at most half of them lie in $A$ for the same
reason. We obtain
\begin{eqnarray*}
n_+(y,A) &\le&|A|\wedge \biggl\lceil\frac{p(y)+1}{2} \biggr\rceil +\sum
_{j=0}^{p(y)-1}\frac{{i_{j+1}}-{i_j}-1}{2}+
\frac{n-i_{p(y)}}{2}
\\
&=&|A|\wedge \biggl\lceil\frac{p(y)+1}{2} \biggr\rceil+\frac{n-p(y)}{2}
\end{eqnarray*}
as required.
\end{pf}
%


\section{A point processes approach}
\label{s_ppp}

In this section, we use point processes techniques to understand the
joint scaling limit of
the random variables $Z_t^{\ssup{1, c}}$, $Z_t^{\ssup{2, c}}$,
$\varPsi
_{t,c}(Z_t^{\ssup{1, c}})$,
$\varPsi_{t,c}(Z_t^{\ssup{2, c}})$ for each $c$ and, in particular,
that of
$Z_t^{\ssup1}$, $Z_t^{\ssup2}$, $\varPsi_t(Z_t^{\ssup1})$, $\varPsi
_t(Z_t^{\ssup2})$.
We show that $Z_t^{\ssup{1, c}}$ and $Z_t^{\ssup{2, c}}$ grow at scale
$r_t$ and that $\varPsi_{t,c}(Z_t^{\ssup{1, c}})-a_{r_t}$ and
$\varPsi_{t,c}(Z_t^{\ssup{2, c}})-a_{r_t}$
grow or decay at scale $d_t$
(which goes to infinity for $\gamma<1$, is a constant for $\gamma=1$,
and tends to zero for $\gamma>1$), and
we find their joint scaling limit in Proposition~\ref{l_4den}. In
particular, we show that the probability
of the event $\mathcal{E}_c(t)$ defined in \eqref{eee} tends to one for
any $c$ and so it suffices to prove complete
localisation and ageing on the event $\mathcal{E}_c(t)$ for a
sufficiently large constant $c$. This constant will be identified
later in Proposition~\ref{Hh} in Section~\ref{s_neg}. Finally, in the
end of this section we prove
Theorem~\ref{main_z}.

For all $z\in\Z^d$ and all sufficiently large $r$, denote
\[
X_{r,z}=\frac{\xi(z)-a_r}{d_r}
\]
and define
\[
\varSigma_r=\sum_{z\in\Z^d}
\varepsilon_{(zr^{-1},X_{r,z})},
\]
where $\varepsilon_x$ denotes the Dirac measure in $x$.
For each $\tau\in\R$ and $q>0$, let
\[
H_{\tau}^q= \bigl\{(x,y)\in\dot{\R}^d\times(-
\infty,\infty]\dvtx y\ge q|x|+\tau \bigr\},
\]
where $\dot{\R}^d$ denotes the one-point compactification of the
Euclidean space. It was proved in \cite{HMS08}, Lemma~4.3,
that for $0<\gamma\le1$ the restriction of each
$\varSigma_r$ to $H_{\tau}^q$ is a point process and, as $r\to\infty$,
$\varSigma_r|_{H_{\tau}^q}$
converges in law to a Poisson point process $\varSigma$ on $H_{\tau}^q$
with intensity measure
\[
\eta(\dd x,\dd y)=\dd x\otimes\gamma e^{-\gamma y}\,\dd y.
\]
However, it is easy to check that the same proof works for all $\gamma>0$.

Observe that we need to restrict $\varSigma_r$ from $\R^d\times\R$ to
$H_{\tau}^q$ in order to ensure
 that there are only finitely many points of $\varSigma_r$ in every
relatively compact set. This is achieved with the help of $q$,
and $\tau$ makes it possible for the spaces $H_{\tau}^q$ to capture the
behaviour of $\varSigma_r$
on the whole space $\R^d\times\R$ as it can be chosen arbitrarily small.

For each $\tau\in\R$ and $\alpha> -\theta$, let
\[
\hat H_{\tau}^{\alpha}= \bigl\{(x,y)\in\dot{\R}^{d+1}
\dvtx y\ge\alpha|x|+\tau \bigr\},
\]
where the hat over $H$ reflects the fact that
the spaces $\dot{\R}^d\times(-\infty,\infty]$ and $\dot{\R}^{d+1}$
have different topology.

For all $c\in\R$, $z\in\Z^d$, and all sufficiently large $t$ define
\[
Y_{t,z,c}=\frac{\varPsi_{t,c}(z)-a_{r_t}}{d_{r_t}} \quad\mbox {and}\quad\varPi_{t,c}=
\sum_{z\in\Z^d}\varepsilon_{(zr_t^{-1},Y_{t,z,c})}.
\]
Recall the definitions of $Y_{t,z}$ and $\varPi_t$ from \eqref{yy} and
\eqref{pp} and observe that $Y_{t,z,c}=Y_{t,z,0}$
and $\varPi_t=\varPi_{t,0}$.

%
\begin{lemma}
\label{l_ppp}
Let $c\in\R$.
For all sufficiently large $t$, $\varPi_{t,c}$ is a point process on~$\hat
H_{\tau}^{\alpha}$.
As $t\to\infty$, $\varPi_{t,c}$ converges in law to a Poisson point
process $\varPi$ on $\hat H_{\tau}^{\alpha}$ with intensity measure
\[
\nu(\dd x,\dd y) =\dd x\otimes\gamma\exp \bigl\{-\gamma\bigl(y+\theta |x|\bigr) \bigr\} \,
\dd y.
\]
\end{lemma}

\begin{pf} Observe that
\[
Y_{t,z,c}=\frac{\xi(z)-a_{r_t}}{d_{r_t}}-\frac{|z|}{\gamma t
d_{r_t}}\log\log t+
\frac{c|z|}{td_{r_t}} =\frac{\xi(z)-a_{r_t}}{d_{r_t}}- \bigl(\theta +o(1) \bigr)
\frac{|z|}{r_t}.
\]
Choose $\alpha'$ and $q$ so that $-\theta<\alpha'<\alpha$ and
$\alpha
'+\theta<q<\alpha+\theta$. Then
%
%
\begin{equation}
\label{pppref} \varPi_{t,c}|_{\hat H^{\alpha}_{\tau}}= \bigl(
\varSigma_{r_t} |_{H_{\tau
}^q}\circ T_{t,c}^{-1}
\bigr) \big|_{\hat H^{\alpha}_{\tau}},
\end{equation}
where $T_{t,c}\dvtx H_{\tau}^q\to\hat H_{\tau}^{\alpha'}$ is such that\vspace*{-1pt}
\[
T_{t,c}\dvtx(x,y)\mapsto\cases{ %
\bigl(x,y-
\bigl(\theta+o(1) \bigr)|x| \bigr), &\quad$\mbox{if }x\neq\infty \mbox{ and }y\neq
\infty$,
\cr
\infty, &\quad$\mbox{otherwise}$. 
}
\]
We define $T\dvtx H_{\tau}^q\to\hat H_{\tau}^{\alpha'}$ by\vspace*{-1pt}
\[
T\dvtx(x,y)\mapsto\cases{ %
\bigl(x,y-\theta|x|\bigr), & \quad$
\mbox{if }x\neq\infty\mbox{ and }y\neq\infty$,
\cr
\infty, &\quad$
\mbox{otherwise}$. 
}
\]
It was proved in \cite{HMS08}, Lemma~2.5, that one can pass to the limit
in \eqref{pppref} as $t\to\infty$ simultaneously
in the mapping $T_{t,c}$ and the point process $\varSigma_{r_t}$ to get\vspace*{-1pt}
\[
\varPi_{t,c}|_{\hat H^{\alpha}_{\tau}}\Longrightarrow \bigl(
\varSigma|_{H_{\tau}^q}\circ T^{-1} \bigr) \big|_{\hat H^{\alpha}_{\tau}}.
\]
Observe that the conditions of that lemma are satisfied as $T$ is
continuous, $H_{\tau}^q$ is compact,
$T_{t,c}\to T$ uniformly on $\{(x,y)\in H_{\tau}^q\dvtx |x|\ge n\}$ as
$t\to
\infty$ for each $n\in\N$, and\vspace*{-1pt}
\[
\eta \bigl\{(x,y)\in H_{\tau}^q\dvtx|x|\ge n \bigr\}\to0
\qquad\mbox{as }n\to\infty
\]
since $\eta(H_{\tau}^q)$ is finite.
Finally, it remains to notice
that $ (\varSigma|_{H_{\tau}^q}\circ T^{-1} ) |_{\hat
H^{\alpha}_{\tau}}$ is a Poisson process with
intensity measure $\eta\circ T^{-1}=\nu$ restricted on $\hat
H^{\alpha
}_{\tau}$.
\end{pf}
%

%
%
\begin{prop}
\label{l_4den}
Let $c\in\R$.
\begin{enumerate}[(b)]
\item[(a)]As $t\to\infty$,\vspace*{-1pt}
\begin{eqnarray*}
&&\biggl(\frac{Z_t^{\ssup{1,c}}}{r_t}, \frac{\varPsi
_{t,c}(Z_t^{\ssup{1,c}})-a_{r_t}}{d_{r_t}}, \frac{Z_t^{\ssup{2,c}}}{r_t},
\frac{\varPsi_{t,c}(Z_t^{\ssup{2,c}})-a_{r_t}}{d_{r_t}} \biggr)
\\
&&\qquad\Longrightarrow \bigl(X^{\ssup1},
Y^{\ssup1}, X^{\ssup2}, Y^{\ssup2} \bigr),
\end{eqnarray*}
where the limit random variable has density\vspace*{-1pt}
\begin{eqnarray*}
&&p(x_1,y_1, x_2, y_2)
\\
&&\qquad=\gamma^2\exp \bigl\{-\gamma\bigl(y_1+y_2+
\theta|x_1|+\theta|x_2| \bigr)-2^d(\gamma
\theta)^{-d} e^{-\gamma y_2} \bigr\}\id_{\{y_1>y_2\}}.
\end{eqnarray*}

\item[(b)]$\Prob\{\mathcal{E}_c(t)\}\to1$ as $t\to\infty$.
\end{enumerate}
\end{prop}

\begin{pf}
(a) Let $A\subset\hat H_{\tau}^0\times\hat H_{\tau}^0$ for some
$\tau$, and
assume that $\operatorname{Leb}(\partial A)=0$. Since $H_{\tau}^0$ is compact,
we have by Lemma~\ref{l_ppp}
%
%
\begin{eqnarray}
\label{4den} &&\Prob \biggl\{ \biggl(\frac{Z_t^{\ssup{1,c}}}{r_t}, \frac{\varPsi_{t,c}(Z_t^{\ssup{1,c}})-a_{r_t}}{d_{r_t}},
\frac{Z_t^{\ssup{2,c}}}{r_t}, \frac{\varPsi_{t,c}(Z_t^{\ssup
{2,c}})-a_{r_t}}{d_{r_t}} \biggr)\in A \biggr\}
\nonumber
\\
&&\qquad=\int_A \id_{\{y_1>y_2\}}\Prob\bigl\{
\varPi_{t,c}(\dd x_1\times\dd y_1)=
\varPi_{t,c}(\dd x_2\times\dd y_2)=1,
\nonumber
\\
&&\phantom{\qquad=\int_A \id_{\{y_1>y_2\}}\Prob \bigl
\{} \varPi_{t,c} \bigl(\R^d\times(y_1,\infty)
\bigr)=\varPi_{t,c} \bigl(\R^d\times(y_2,y_1)
\bigr)=0 \bigr\}
\\
&&\qquad\to\int_A \id_{\{y_1>y_2\}}\Prob \bigl\{
\varPi(\dd x_1\times\dd y_1)=1 \bigr\}\Prob \bigl\{\varPi(
\dd x_2\times\dd y_2)=1 \bigr\}
\nonumber
\\
&&\phantom{\qquad\to\int_A} {} \times\Prob \bigl\{\varPi
\bigl(\R^d\times(y_1,\infty) \bigr)=0 \bigr\}\Prob \bigl\{
\varPi \bigl(\R^d\times(y_2,y_1) \bigr)=0
\bigr\}
\nonumber
\\
&&\qquad=\int_A \id_{\{y_1>y_2\}}\nu(\dd
x_1,\dd y_1)\nu(\dd x_2,\dd y_2)
\exp \bigl\{-\nu \bigl(\R^d\times(y_2,\infty) \bigr) \bigr
\}.
\nonumber
\end{eqnarray}
Integrating we obtain
%
%
\begin{eqnarray}
\label{dom} \nu \bigl(\R^d\times(y_2,\infty) \bigr)
&=& \gamma\int_{\R^d}\int_{y_2}^{\infty}
\exp\bigl \{-\gamma y-\gamma\theta|x|\bigr \} \,\dd y\,\dd x
\nonumber
\\[-8pt]
\\[-8pt]
&=&2^d(\gamma\theta)^{-d}e^{-\gamma y_2}.
\nonumber
\end{eqnarray}
Substituting this, as well as the expressions for $\nu(\dd x_1,\dd
y_1)$ and $\nu(\dd x_2,\dd y_2)$ into \eqref{4den}
we obtain
\begin{eqnarray*}
&&\lim_{t\to\infty}\Prob \biggl\{ \biggl(\frac{Z_t^{\ssup{1,c}}}{r_t},
\frac{\varPsi_{t,c}(Z_t^{\ssup{1,c}})-a_{r_t}}{d_{r_t}}, \frac
{Z_t^{\ssup{2,c}}}{r_t}, \frac{\varPsi_{t,c}(Z_{t,c}^{\ssup
{2,c}})-a_{r_t}}{d_{r_t}} \biggr\}\in A \biggr)
\\
&&\qquad= \int_A p(x_1,y_1,x_2,y_2)
\,\dd x_1\,\dd y_1\,\dd x_2\,\dd
y_2.
\end{eqnarray*}

It remains now to generalise this equality to all sets $A\subset\R
^d\times\R$ with $\operatorname{Leb}(\partial A)=0$.
Since $\tau$ can be arbitrarily small, to do so it suffices to show
that $p$ integrates to one. We have
%
%
\begin{eqnarray}
\label{integrate} &&\int_{\R^d\times\R\times\R^d\times\R
}p(x_1,x_2,y_1,y_2)
\,\dd x_1\,\dd y_1\,\dd x_2\,\dd
y_2
\nonumber
\\
&&\qquad=2^{2d}(\gamma\theta)^{-2d}\int_{-\infty}^\infty
\int_{y_2}^{\infty} \gamma^2\exp \bigl\{-
\gamma(y_1+y_2)-2^d (\gamma\theta
)^{-d}e^{-\gamma
y_2} \bigr\}\,\dd y_1\,\dd
y_2
\nonumber
\\[-8pt]
\\[-8pt]
&&\qquad=2^{2d}(\gamma\theta)^{-2d}\int_{-\infty}^\infty
\gamma\exp \bigl\{-2\gamma y_2-2^d (\gamma
\theta)^{-d}e^{-\gamma
y_2} \bigr\}\,\dd y_2
\nonumber
\\
&&\qquad=\int_{0}^\infty ue^{-u}\,\dd
u=1,
\nonumber
\end{eqnarray}
where in the last line we used the substitution $u=2^d (\gamma\theta
)^{-d}e^{-\gamma y_2}$.

(b) This immediately follows from (a) since $d_{r_t}=d_t (1+o(1))$ and
$f_t\to0$, $g_t\to\infty$, $\lambda_t\to0$.
\end{pf}
%

%
\begin{remark}
\label{rem3}
The reason why we need to study a general $c$ rather than $c=0$ is just
to show that
$|Z^{\ssup{1,c}}_t|<r_tg_t$ and $|Z^{\ssup{2,c}}_t|<r_tg_t$ with high
probability, which is done in part (b) of the proposition above.
This will be required
later on in Lemma~\ref{l_u2} with some $c$ identified in Proposition~\ref{Hh}. The full strength of the convergence result
proved in the part (a) of the proposition will only be used for $c=0$.
\end{remark}

\begin{pf*}{Proof of Theorem~\ref{main_z}}
The result follows from Proposition~\ref{l_4den}(a) with $c=0$ by
integrating the density $p$ over all possible values
of $x_2$, $y_1$, and $y_2$. Similarly to \eqref{integrate}, we obtain
\begin{eqnarray*}
p^{\ssup1}(x) & =&\int_{\R\times\R^d\times\R}p(x,y_1,x_2,y_2)
\,\dd y_1\,\dd x_2\,\dd y_2
\\
&=&2^{d}(\gamma\theta)^{-d}\exp\bigl \{-\gamma\theta|x| \bigr \}
\\
&&{} \times\int_{-\infty}^\infty\int
_{y_2}^{\infty} \gamma^2\exp \bigl\{-
\gamma(y_1+y_2)-2^d (\gamma\theta
)^{-d}e^{-\gamma
y_2} \bigr\}\,\dd y_1\,\dd
y_2
\\
&=&2^{-d}d^{d(1-1/\gamma)}\exp \bigl\{-d^{1-1/\gamma}|x| \bigr\}
\end{eqnarray*}
as required.
\end{pf*}



\section{Negligible paths of the random walk}
\label{s_neg}

Throughout this section, we assume that $0<\gamma<2$.
We introduce three groups of paths of the random walk $(X_s)$
informally described in the \hyperref[intro]{Introduction} and show that
their contribution to the total mass of the solution $u$ of the
parabolic Anderson model is negligible.

Denote by $J_t$ the number of jumps the random walk $(X_s)$ makes up to
time $t$ and
consider the following three groups of paths:
\[
E_i(t)=\cases{ %
\Bigl\{\displaystyle\max
_{0\leq s\leq t}{\xi(X_s)}=\xi\bigl(Z_t^{\ssup1}
\bigr), J_t > \bigl |Z_t^{\ssup
1}\bigr |(1+\rho_t)
\Bigr\}, &\quad$i=1$,\vspace*{2pt}
\cr
\Bigl\{\xi_{r_tg_t}^{\ssup{k_t}} \le\displaystyle\max
_{0\leq s\leq t}{\xi (X_s)}\neq\xi\bigl(Z_t^{\ssup1}
\bigr) \Bigr\}, & \quad$i=2$,\vspace*{2pt}
\cr
\Bigl\{\displaystyle\max_{0\leq s\leq t}{
\xi(X_s)}< \xi_{r_tg_t}^{\ssup{k_t}} \Bigr\} , & \quad$i=3$.
}
\]
Denote by
\[
U_i(t)=\E_0 \biggl[\exp \biggl\{\int
_0^t\xi(X_s)\,\dd s \biggr\}\1
_{E_i(t)} \biggr],\qquad1\le i\le3
\]
their contributions to the total mass of the solution.
The aim of this section is to show that all $U_i(t)$ is negligible with
respect to $U(t)$.

We start with Lemma~\ref{app} where we collect all asymptotic
properties of the environment
which we use later on. In Lemma~\ref{l_lb}, we prove a simple lower
bound for the total mass
$U(t)$. Then we prove Proposition~\ref{Hh}, which is a crucial tool for
analysing $U_1(t)$ and $U_2(t)$
as it gives a general upper bound on the total mass corresponding to
the paths reaching the maximum
of the potential in a certain set
and having a lower bound restriction on the number of jumps $J_t$.
Equipped with this result,
we show that $U_1(t)$ and $U_2(t)$ are negligible in Lemmas~\ref{l_u1}
and~\ref{l_u2}.
Finally, Lemma~\ref{l_u3} provides a simple proof of the negligibility
of $U_3(t)$.

Observe that Proposition~\ref{Hh} identifies the constant $c$, which is
then fixed and used throughout the
paper afterward.\vadjust{\goodbreak}

%
\begin{lemma}
\label{app}
Almost surely,
\begin{enumerate}[(b)]
\item[(a)] $\xi_r^{\ssup{\lfloor r^{\rho}\rfloor}}\sim((d-\rho
)\log
r)^{1/\gamma}$
and $ \xi_r^{\ssup{\lfloor r^{\sigma}\rfloor}}\sim((d-\sigma)\log
r)^{1/\gamma}$ as $r\to\infty$;
\item[(b)] $\xi_{r_tg_t}^{\ssup{k_t}}\sim((d-\rho)\log
t)^{1/\gamma
}$ and
$ \xi_{r_tg_t}^{\ssup{m_t}}\sim((d-\sigma)\log t)^{1/\gamma}$ as
$t\to
\infty$;
\item[(c)] $\log(\xi_r^{\ssup{\lfloor r^{\rho}\rfloor}}-\xi
_r^{\ssup
{\lfloor r^{\sigma}\rfloor}})= \frac{1}{\gamma}\log\log r+O(1)$
as $r\to\infty$;
\item[(d)] $\log(\xi_{r_tg_t}^{\ssup1}-\xi_{r_tg_t}^{\ssup
{m_t}})=\frac{1}{\gamma}\log\log t+O(1)$ as $t\to\infty$;
\item[(e)] the set $G_p$ is totally disconnected eventually for all $p$.
\end{enumerate}
Further,
\begin{enumerate}[(g)]
\item[(f)] for all $c$, $Z_t^{\ssup1}\in F_{r_tg_t}$ on the event
$\mathcal
{E}_c(t)$ eventually for all $t$;
\item[(g)] for all $c$, $\log\xi(Z_t^{\ssup1})=\frac{1}{\gamma
}\log\log
t+O(1)$ on the event $\mathcal{E}_c(t)$ as $t\to\infty$;
\item[(h)] there exists a constant $c_1>0$ such that
$|z|>t^{c_1}$ for all $z\in F_{r_tg_t}$ eventually for all $t$ almost surely.
\end{enumerate}
%
%
\end{lemma}

\begin{pf} (a) It follows from the proof of \cite{HMS08}, Lemma~4.7, that for each $\kappa\in(0,d)$ almost surely
\[
\xi_r^{\ssup{\lfloor r^{\kappa}\rfloor}}\sim \bigl((d-\kappa)\log r
\bigr)^{1/\gamma}
\]
as $r\to\infty$. It remains to substitute $\kappa=\rho$ and $\kappa
=\sigma$.

(b) This follows from (a) since $k_t=\lfloor(r_tg_t)^{\rho}\rfloor$
and $m_t=\lfloor(r_tg_t)^{\sigma}\rfloor$.

(c) This follows from (a) since $\rho\neq\sigma$.

(d) This follows from (a) and Lemma~\ref{l_asasas} since $\rho\neq0$.

(e) This was proved in \cite{KLMS09}, Lemma~2.2, for Pareto potentials
(observe that the proof relies on
$\sigma<1/2$ which is the reason why we have imposed this restriction).
It remains to notice that
$\xi(z)=(\alpha\log(\zeta(z)))^{1/\gamma}$, where $\{\zeta
(z)\dvtx z\in\Z
^d\}$
is a Pareto-distributed potential with parameter $\alpha$. As the
locations of upper order statistics for $\zeta$ and $\xi$ coincide,
we obtain that $G_p$ is eventually totally disconnected for Weibull
potentials as well.

(f) Denote by $w_t$ the maximiser of $\xi$ in the ball of radius $t$.
Using Lemma~\ref{l_asasas}, we obtain
\begin{eqnarray*}
\xi\bigl(Z_t^{\ssup1}\bigr) &\ge&\varPsi_t
\bigl(Z_t^{\ssup1}\bigr)\ge\varPsi_t(w_t)
=\xi(w_t)-\frac{|w_t|}{\gamma t}\log\log t
\\
&\ge&\xi_t^{\ssup1}-\frac{1}{\gamma}\log\log t \sim(d\log
t)^{1/\gamma}.
\end{eqnarray*}
It remains to observe that
$|Z_t^{\ssup1}|\le r_tg_t$ on the event $\mathcal{E}_c(t)$ and use
(a) to get
\[
\xi\bigl(Z_t^{\ssup1}\bigr)\ge \bigl((d-\rho)\log t
\bigr)^{1/\gamma}\sim \xi_{r_tg_t}^{\ssup{k_t}}.
\]

(g) It follows from (f) that $\log\xi_{r_tg_t}^{\ssup{k_t}}\le
\log
\xi(Z_t^{\ssup1})\le\log\xi_{r_tg_t}^{\ssup1}$ on the
event $\mathcal{E}_c(t)$. It remains to
observe that
$\log\xi_{r_tg_t}^{\ssup{k_t}}=\frac{1}{\gamma}\log\log t+O(1)$
according to (a) and
$\log\xi_{r_tg_t}^{\ssup1}=\frac{1}{\gamma}\log\log t+O(1)$ by Lemma~\ref{l_asasas}.

(h)
Choose $c_1$ small enough so that $c_1(d+c_1)<d-\rho-c_1$.
Then almost surely eventually
\[
\xi^{\ssup1}_{t^{c_1}}\le \bigl((d+c_1)\log
t^{c_1} \bigr)^{1/\gamma}< \bigl((d-\rho-c_1)\log t
\bigr)^{1/\gamma}<\xi_{r_tg_t}^{\ssup{k_t}},
\]
which implies the result.
\end{pf}
%

%
\begin{lemma}
\label{l_lb}
For each $c$,
%
%
\begin{equation}
\label{lb} \log U(t) \ge t\varPsi_t\bigl(Z_t^{\ssup1}
\bigr)-2dt+O(r_tg_t)
\end{equation}
on the event $\mathcal{E}_c(t)$
eventually for all $t$.
\end{lemma}

\begin{pf}
The idea of the proof is the same as of \cite{HMS08}, Lemma~2.1, for
Weibull potentials
and \cite{KLMS09}, Proposition~4.2, for Pareto potentials.
However, we need to estimate the error term more precisely.

Let $\rho\in(0,1]$ and $z\in\Z^d$, $z\neq0$. Following the lines of
\cite{KLMS09}, Proposition~4.2,
we obtain
%
%
\begin{equation}
\label{l1} U(t)\ge\exp \biggl\{t(1-\rho)\xi(z)-|z|\log\frac
{|z|}{e\rho
t}-2dt+O\bigl(
\log|z|\bigr) \biggr\}.
\end{equation}
Take $z=Z_t^{\ssup1}$ and $\rho=|Z_t^{\ssup1}|/(t\xi(Z_t^{\ssup
1}))$. Observe that on the event
$\mathcal{E}_c(t)$ this $\rho$ belongs to $(0,1]$
eventually as
\[
\frac{|Z_t^{\ssup1}|}{t\xi(Z_t^{\ssup1})}\le\frac{r_tg_t}{t\xi
_{r_tg_t}^{\ssup{k_t}}} =O \biggl(\frac{g_t}{\log t\cdot\log\log t}
\biggr)=o(1)
\]
by Lemma~\ref{app}(f) and according to \eqref{a}. Substituting this
into \eqref{l1} and using Lemma~\ref{app}(g) we obtain
\begin{eqnarray*}
\log U(t) &\ge& t\xi\bigl(Z_t^{\ssup1}
\bigr)-\bigl |Z_t^{\ssup1}\bigr |\log\xi \bigl(Z_t^{\ssup1}
\bigr)-2dt+O(\log t)
\\
&=&t\varPsi_t\bigl(Z_t^{\ssup1}
\bigr)-2dt+O(r_tg_t)
\end{eqnarray*}
on the event $\mathcal{E}_c(t)$.
\end{pf}

For all sufficiently large $t$, consider a set $M_t\subset\Z^d$ and
a nonnegative function $h_t=O(r_tg_t)$ (which may both depend on $\xi$).
Denote by $z_t$ a point along the trajectory of $(X)_s$, $s\in[0,t]$,
where the value of the potential is maximal.
Define
\[
U_{M,h}(t)=\E_0 \biggl[\exp \biggl\{\int
_0^t\xi(X_s)\,\dd s \biggr\} \id
\Bigl\{ \max_{0\leq s\leq t}{\xi(X_s)}\ge
\xi_{r_tg_t}^{\ssup{k_t}}, z_t\in M_t,
J_t\ge h_t \Bigr\} \biggr].
\]
%
In the sequel, $U_{M,h}(t)$ will correspond to $U_1(t)$ if we choose
$M_t=\{Z_t^{\ssup1}\}$, ${h_t=|Z_t^{\ssup1}|(1+\rho_t)}$ and to
$U_2(t)$ if we
choose $M_t=\Z^d\backslash\{Z_t^{\ssup1}\}$, $h_t=0$.

%
%
\begin{prop}
\label{Hh} There is a constant $c$ such that
\begin{eqnarray*}
&&\log U_{M,h}(t) \le\max\biggl\{t\varPsi_t\bigl(Z_t^{\ssup2}
\bigr),
\\
&&\phantom{\log U_{M,h}(t) \le\max \biggl\{} {}\max_{z\in M_t} \biggl\{ t\varPsi_{t,c}(z)
-\frac{(h_t-|z|)_+}{2}
\bigl(\gamma^{-1}-\beta \bigr)\log\log t \biggr\}
\\
&&\phantom{\hspace*{257pt}}{}+O(r_tg_t)
\biggr\}
\\
&&\phantom{\log U_{M,h}(t) \le}{}-2dt
\end{eqnarray*}
on the event $\mathcal{E}_c(t)$ eventually for all $t$.
\end{prop}

\begin{pf}
Consider the event $\mathcal{E}_c(t)$ and suppose that $t$ is
sufficiently large.
Using the notation from Section~\ref{s_geo},
for each $n, p\in\N\cup\{0\}$ and $t$ large enough, we denote
\[
\mathcal{P}_{n,p}(t)= \bigl\{y\in\mathcal{P}_n\dvtx
y_0=0, p(y)=p, q(y)>\xi_{r_tg_t}^{\ssup{k_t}}, z(y)\in
M_t \bigr\}.
\]
Observe that
$q(y)\ge\xi_{r_tg_t}^{\ssup{k_t}}$ implies by Lemma~\ref{app}(h) that
$p(y)>t^{c_1}$,
for some \mbox{$c_1>0$}. In particular,
%
%
\begin{equation}
\label{logp} \log\log p(y)\ge\log\log t+\log c_1.
\end{equation}
We have
\[
U_{M,h}(t)=\sum_{n\ge h_t}\sum
_{t^{c_1}<p\le n}\sum_{y\in\mathcal
{P}_{n,p}(t)}U(t,y),
\]
where $U(t,y)$ has been defined in \eqref{uty}.
Since the number of paths in the set $\mathcal{P}_{n,p}(t)$ is bounded
by $(2d)^n$, we obtain
\begin{eqnarray*}
U_{M,h}(t) &\le&\sum_{p>t^{c_1}}\sum
_{n\ge p\vee h_t}(2d)^{-n}\max_{y\in
\mathcal
{P}_{n,p}(t)} \bigl
\{(2d)^{2n}U(t,y) \bigr\}
\\
&\le&4\max_{p>t^{c_1}}\max_{n\ge p\vee h_t} \max
_{y\in\mathcal
{P}_{n,z}(t)} \bigl\{(2d)^{2n}U(t,y) \bigr\}
\end{eqnarray*}
and so
%
%
\begin{equation}
\label{prep} \log U_{M,h}(t) \le\max_{p>t^{c_1}}\max
_{n\ge p\vee h_t} \max_{y\in\mathcal{P}_{n,z}(t)} \bigl\{3n\log (2d)+\log
U(t,y) \bigr\}.
\end{equation}
Let $p>t^{c_1}$, $n\ge p\vee h_t$, and $y\in\mathcal{P}_{n,p}(t)$.
Denote $i(y)=\min\{i\dvtx\xi(y_i)=q(y)\}$
and
%
%
\begin{equation}
\label{o} Q(p,y)=q(y)\vee\xi_p^{\ssup{\lfloor p^{\rho}\rfloor}}+
\mu_p,
\end{equation}
where the correction term $\mu_p$ has been defined in \eqref{mu}.
Define
\[
\xi^{y}_i=\cases{ %
\xi(y_i),&\quad$\mbox{if }i\neq i(y)$,
\cr
Q(p,y),&\quad$\mbox{if }i=
i(y)$. 
}
\]
Since $\xi^y_i\ge\xi(y_i)$ for all $i$, we have
\begin{eqnarray*}
&&U(t,y)\le(2d)^{-n}\mathsf{E} \Biggl[\exp \Biggl\{\sum
_{i=0}^{n-1}\tau_i\xi^y_i+
\Biggl(t-\sum_{i=0}^{n-1}\tau_i
\Biggr)\xi^y_n \Biggr\}
\\
&&\phantom{\hspace*{112pt}} {}\times\one \Biggl\{\sum
_{i=0}^{n-1}\tau_i<t,\sum
_{i=0}^{n}\tau_i>t \Biggr\} \Biggr].
\end{eqnarray*}
This expectation has been bounded from above in (4.16) and (4.17) of
\cite{HMS08}. Substituting its bound, we obtain
\begin{eqnarray*}
U(t,y) &\le&\exp \bigl\{t\xi_{i(y)}^{y}-2dt \bigr\}\prod
_{i\neq i(y)}\frac
{1}{\xi
_{i(y)}^{y}-\xi_i^y}
\\
&=&\exp \bigl\{tQ(p,y)-2dt \bigr\}\prod_{i\neq i(y)}
\frac{1}{Q(p,y)-\xi(y_i)}
\end{eqnarray*}
and hence
%
%
\begin{equation}
\label{e} \log U(t,y)\le tQ(p,y)-2dt-\sum_{i\neq i(y)}
\log \bigl(Q(p,y)-\xi(y_i) \bigr).
\end{equation}

The set $G_{p}$ consists of $\lfloor p^{\sigma}\rfloor$
elements and is totally disconnected by Lem\-ma~\ref{app}(e). Hence, by
Lemma~\ref{l_count} we have
%
%
\begin{equation}
\label{o1} n_{+}(y,G_{p})\le\frac{n-p}{2}+p^{\sigma}.
\end{equation}
In each point $y_i\in G_p$ we use \eqref{o} to estimate
%
%
\begin{equation}
\label{o2} \log \bigl(Q(p,y)-\xi(y_i) \bigr) \ge\log
\mu_p= -\beta\log\log p.
\end{equation}
On the other hand,
%
%
\begin{eqnarray}
\label{o3} n_{-}(y,G_{p})&=&n+1-n_{+}(y,G_{p})
\nonumber
\\[-8pt]
\\[-8pt]
&\ge& n+1-\frac{n-p}{2}-p^{\sigma} =p-p^{\sigma}+
\frac{n-p}{2}+1\nonumber
\end{eqnarray}
and in each point $y_i\notin G_p$ we obtain by Lemma~\ref{app}(c)
%
%
\begin{equation}
\label{o4} \log \bigl(Q(p,y)-\xi(y_i) \bigr)\ge\log \bigl(
\xi_p^{\ssup{\lfloor
p^{\rho
}\rfloor}}-\xi_p^{\ssup{\lfloor p^{\sigma}\rfloor}} \bigr) \ge
\gamma^{-1}\log\log p+c_2
\end{equation}
with some constant $c_2$.
Using \eqref{o2} and \eqref{o4} together with \eqref{e},
we obtain
\begin{eqnarray*}
\log U(t,y) &\le& tQ(p,y)-2dt +n_{+}(y,G_{p})\beta\log
\log p
\\
&&{}- \bigl(n_{-}(y,G_p)-1 \bigr) \bigl(
\gamma^{-1}\log\log p+c_2 \bigr).
\end{eqnarray*}
Substituting \eqref{o1} and \eqref{o3} and using $p^{\sigma}\log
\log
p\le n$, we obtain
%
%
\begin{eqnarray}
\label{uu5} && 3n\log(2d)+\log U(t,y)
\nonumber
\\
&&\qquad\le3n\log(2d)+tQ(p,y)-2dt
\nonumber
\\[-8pt]
\\[-8pt]
&&\qquad\phantom{\le} {} + \biggl[\frac{n-p}{2}+p^{\sigma} \biggr]
\beta\log\log p - \biggl[p-p^{\sigma}+\frac{n-p}{2} \biggr] \bigl(
\gamma^{-1}\log\log p+c_2 \bigr)
\nonumber
\\
&&\qquad\le tQ(p,y)-\frac{p}{\gamma}\log\log p-2dt -\frac{n-p}{2} \bigl(
\gamma^{-1}-\beta \bigr)\log\log p+c_3n
\nonumber
\end{eqnarray}
with some constant $c_3$.

Now we distinguish between the following two cases.

\textit{Case 1}. Suppose $q(y)\ge\xi_p^{\ssup{\lfloor p^{\rho
}\rfloor
}}$. Then $Q(p,y)=\xi(z(y))+\mu_p$
and estimating $p\ge|z(y)|$ we get
\begin{eqnarray*}
3n\log(2d)+\log U(t,y) &\le& t\xi \bigl(z(y) \bigr)+t
\mu_p-\frac{|z(y)|}{\gamma}\log\log p-2dt
\\
&&{}-\frac{n-|z(y)|}{2} \bigl(\gamma^{-1}-\beta \bigr)\log\log p
+c_3n.
\end{eqnarray*}
Observe that $t\mu_p\le t\mu_{t^{c_1}}= t(c_1\log t)^{-\beta
}=o(r_tg_t)$ since $\beta>1-1/\gamma$
and according to \eqref{a}.
Using monotonicity in $n$ and $n\ge|z(y)|\vee h_t$ together with
\eqref
{logp}, we obtain
%
%
\begin{eqnarray}
\label{u7} &&3n\log(2d)+\log U(t,y)
\nonumber
\\
&&\qquad\le t\varPsi_t \bigl(z(y) \bigr)+c\bigl |z(y)\bigr |-2dt
\nonumber
\\[-8pt]
\\[-8pt]
&&\qquad\phantom{\le} {}-\frac{(h_t-|z(y)|)_+}{2} \bigl(\gamma ^{-1}-\beta
\bigr)\log\log t +c h_t+o(r_tg_t)
\nonumber
\\
&&\qquad\le\max_{z\in M_t} \biggl\{ t\varPsi_{t,c}(z)-
\frac
{(h_t-|z|)_+}{2} \bigl(\gamma^{-1}-\beta \bigr)\log\log t \biggr\}-2dt
+O(r_tg_t)
\nonumber
\end{eqnarray}
with some constant $c$.

\textit{Case 2}. Suppose $q(y)< \xi_p^{\ssup{\lfloor p^{\rho
}\rfloor
}}$. Then $Q(p,y)= \xi_p^{\ssup{\lfloor p^{\rho}\rfloor}}+\mu_p$.
Now \eqref{uu5} implies
\begin{eqnarray*}
3n\log(2d)+\log U(t,y) &\le& t\xi_p^{\ssup{\lfloor p^{\rho}\rfloor}}+t
\mu_p-\frac
{p}{\gamma}\log\log p-2dt
\\
&&{}-\frac{n-p}{2} \bigl(\gamma^{-1}-\beta \bigr)\log\log p
+c_4n
\end{eqnarray*}
with some constant $c_4$. Using monotonicity in $n$ and $n\ge p$, we get
\[
3n\log(2d)+\log U(t,y) \le t\xi_p^{\ssup{\lfloor p^{\rho}\rfloor
}}+t(\log
p)^{-\beta
}-\frac
{p}{\gamma}\log\log p-2dt+c_4p.
\]
By Lemma~\ref{app}(a) and using $\beta\ge0$, we obtain that the second
term is dominated by the first one, the fifth by the third one, and so
%
%
\begin{equation}
\label{u6} \quad 3n\log(2d)+\log U(t,y) \le t \bigl((d-\rho/2)\log p
\bigr)^{1/\gamma}-c_5p\log\log p-2dt
\end{equation}
with some constant $c_5>0$. Differentiating, we obtain the following
equation for the maximiser
$p_t$ of the expression on the right-hand side of \eqref{u6}:
\[
\frac{t(d-\rho/2) ((d-\rho/2)\log p_t)^{1/\gamma-1}}{\gamma
p_t}-c_5\log\log p_t-
\frac{c_5}{\log p_t}=0.
\]
Resolving this asymptotics, we obtain
\[
p_t=r_t(d-\rho/2)^{1/\gamma} \bigl(1+o(1) \bigr).
\]
%
Finally, substituting this into \eqref{u6} yields
%
%
\begin{eqnarray}
\label{u8} 3n\log(2d)+\log U(t,y) &\le& t \bigl((d-\rho/3)\log
r_t \bigr)^{1/\gamma}-2dt
\nonumber
\\
&\le& \bigl(1-\rho/(3d) \bigr)^{1/\gamma}ta_{r_t}-2dt
\\
&\le& t\varPsi_t\bigl(Z_t^{\ssup2}\bigr)-2dt
\end{eqnarray}
on the event $\mathcal{E}_c(t)$.
It remains to substitute \eqref{u7} and \eqref{u8} into \eqref{prep} to
complete the proof.
\end{pf}
%

%
\begin{remark}
\label{rem1}
Observe that the scaling function $\mu_p$, being part of $Q(p,y)$,
appears both in the main and in the logarithmic term
of \eqref{e}. Being part of the main term, $t\mu_p$ needs to be as
small as $O(r_tg_t)$ in order to not imbalance the significant terms.
This leads to the restriction $\beta>1-1/\gamma$. However, as a part of
the logarithmic term, $\mu_p$ needs to be large enough
so that the contribution $\gamma^{-1}\log\log p$ of ``good'' points
$y_i\notin G_p$ dominates over the contribution $\beta\log\log p$
of ``bad'' points $y_i\in G_p$. This imposes the restriction $\beta
<1/\gamma$. The combination of these two conditions only allows
to choose such $\beta$ if $0<\gamma<2$.
\end{remark}

From now on, we assume that the constant $c$ is fixed and chosen
according to Proposition~\ref{Hh}.

%
\begin{lemma}
\label{l_u1}
Almost surely,
\[
\frac{U_1(t)}{U(t)}\1_{\mathcal{E}_c(t)}\to0 \qquad\mbox{as }t\to \infty.
\]
\end{lemma}

\begin{pf}
We use Proposition~\ref{Hh} with $M_t=\{Z_t^{\ssup1}\}$ and
$h_t=|Z_t^{\ssup1}|(1+\rho
_t)$. Clearly $h_t=O(r_tg_t)$
on the event $\mathcal{E}_c(t)$. By Lemma~\ref{app}(f),
we have $Z_t^{\ssup1}\in F_{r_tg_t}$, which implies
$U_{M,h}(t)=U_1(t)$ eventually for all $t$. Since $|Z_t^{\ssup1}|\le
r_tg_t$ and
so $t\varPsi_{t,c}(Z_t^{\ssup1})=t\varPsi_t(Z_t^{\ssup
1})+O(r_tg_t)$, we obtain
%
%
\begin{eqnarray}
\label{u10} \log U_1(t) &\le&\max \biggl\{ t\varPsi_t
\bigl( Z_t^{\ssup2}\bigr), t\varPsi_{t}
\bigl(Z_t^{\ssup1}\bigr)-\frac{|Z_t^{\ssup1}|\rho
_t}{2}(1/\gamma-\beta)\log
\log t +O(r_tg_t) \biggr\}
\nonumber
\\[-5pt]
\\[-8pt]
&&{}-2dt.
\nonumber
\end{eqnarray}
In order to show that
%
%
\begin{equation}
\label{u222} \log U_1(t)-\log U(t)\to-\infty
\end{equation}
we consider the terms under the maximum in \eqref{u10} separately.
Using the lower bound for the total mass given by Lemma~\ref{l_lb} and
taking into account that
$\varPsi_t(Z_t^{\ssup1})-\varPsi_t(Z_t^{\ssup2})>d_t\lambda_t$ on
the event $\mathcal{E}_c(t)$,
we get for the first term
%
%
\begin{eqnarray}
\label{u11} t\varPsi_t\bigl(Z_t^{\ssup2}
\bigr)-2dt-\log U(t) &\le& t \varPsi_t\bigl(Z_t^{\ssup2}
\bigr)- t\varPsi_t\bigl(Z_t^{\ssup1}
\bigr)+O(r_tg_t)
\nonumber
\\[-8pt]
\\[-8pt]
&<&-td_t\lambda_t+O(r_tg_t)\to-
\infty
\nonumber
\end{eqnarray}
according to \eqref{a}. For the second term, we again use the lower
bound from Lemma~\ref{l_lb}
and take into account that $|Z_t^{\ssup1}|\ge r_tf_t$ on the event
$\mathcal
{E}_c(t)$. This implies
%
%
\begin{eqnarray}
\label{u12} &&t\varPsi_{t}\bigl(Z_t^{\ssup1}
\bigr)-\frac{|Z_t^{\ssup1}|\rho
_t}{2} (1/\gamma- \beta)\log\log t+O(r_tg_t)-2dt-
\log U(t)
\nonumber
\\
&&\qquad\le-\frac{|Z_t^{\ssup1}|\rho_t}{2}(1/\gamma-\beta)\log \log t+O(r_tg_t)
\\
&&\qquad\le-\frac{r_t f_t \rho_t}{2}(1/\gamma-\beta)\log\log t+O(r_tg_t)
\to-\infty
\nonumber
\end{eqnarray}
by \eqref{a}.
Combining \eqref{u11}, \eqref{u12} and \eqref{u10} we get \eqref{u222}
on the event $\mathcal{E}_c(t)$.
\end{pf}
%

%
\begin{lemma}
\label{l_u2}
Almost surely,
\[
\frac{U_2(t)}{U(t)}\1_{\mathcal{E}_c(t)}\to0 \qquad\mbox{as }t\to \infty.
\]
\end{lemma}

\begin{pf} We use Proposition~\ref{Hh} with $M_t=\Z^d\backslash\{
Z_t^{\ssup1}
\}$ and $h_t=0$. In this case $U_{M,h}(t)=U_2(t)$,
and we have
%
%
\begin{equation}
\label{u20} \log U_2(t)\le\max \Bigl\{t\varPsi_t
\bigl( Z_t^{\ssup2}\bigr), t \max_{z\neq Z_t^{\ssup1}}
\varPsi_{t,c}(z)+O(r_tg_t) \Bigr\}-2dt.
\end{equation}

Since $|Z_t^{\ssup{1, c}}|\le r_tg_t$ and $|Z_t^{\ssup{2, c}}|\le
r_tg_t$ on the event $\mathcal{E}_c(t)$, we have
for $i\in\{1,2\}$
\[
t\varPsi_{t,c} \bigl(Z_t^{\ssup{i, c}} \bigr)=t
\varPsi_t \bigl(Z_t^{\ssup{i,
c}} \bigr)+c|Z_t^{\ssup
{i, c}}|=t
\varPsi_t \bigl(Z_t^{\ssup{i, c}} \bigr)+O(r_tg_t).
\]
Substituting this into \eqref{u20} and observing that $z\neq
Z_t^{\ssup1}$, we obtain
\[
\log U_2(t)\le t\varPsi_t\bigl(Z_t^{\ssup2}
\bigr)+O(r_tg_t)-2dt.
\]
Using the lower bound for the total mass given by Lemma~\ref{l_lb} and
taking into account that
$\varPsi_t(Z_t^{\ssup1})-\varPsi_t(Z_t^{\ssup2})>d_t\lambda_t$ on
the event $\mathcal{E}_c(t)$,
we get
\begin{eqnarray*}
\log U_2(t)-\log U(t) &\le& t\varPsi_t
\bigl(Z_t^{\ssup2}\bigr)-t \varPsi_t
\bigl(Z_t^{\ssup1}\bigr)+O(r_tg_t)
\\
&\le&-t d_t\lambda_t +O(r_tg_t)
\to-\infty
\end{eqnarray*}
according to \eqref{a} on the event $\mathcal{E}_c(t)$.
\end{pf}
%

%
\begin{lemma}
\label{l_u3}
Almost surely,
\[
\frac{U_3(t)}{U(t)}\1_{\mathcal{E}_c(t)}\to0 \qquad\mbox{as }t\to \infty.
\]
\end{lemma}

\begin{pf} We can estimate the integral in the Feynman--Kac formula
for $U_3(t)$ by $t\xi_{r_tg_t}^{\ssup{k_t}}$
and get
\[
\log U_3(t)\le t\xi_{r_tg_t}^{\ssup{k_t}}\sim t \bigl((d-
\rho)\log t \bigr)^{1/\gamma}\le(1-\delta) t a_{r_t}
\]
with some $\delta>0$ eventually for all $t$ by Lemma~\ref{app}(b).
Using the lower bound for $U(t)$ from Lemma~\ref{l_lb}, we have
\begin{eqnarray*}
\log U_3(t)-\log U(t) &\le&(1-\delta)ta_{r_t}-t
\varPsi_t\bigl(Z_t^{\ssup1}\bigr)+2dt+O(r_tg_t)
\\
&\le&-\delta ta_{r_t}+td_t g_t+2dt+O(r_tg_t)
\to-\infty
\end{eqnarray*}
since $\varPsi_t(Z_t^{\ssup1})>a_{r_t}-d_t g_t$ on the event
$\mathcal{E}_c(t)$.
\end{pf}
%


\section{Localisation}
\label{s_loc}

The aim of this section is to prove Theorem~\ref{main_w}.
We assume throughout this section that $0<\gamma<2$ and we suppose
that $c$
is chosen according to Proposition~\ref{Hh}.

Let
\[
B_t= \bigl\{z\in\Z^d\dvtx|z|\le\bigl |Z_t^{\ssup1}\bigr |(1+
\rho_t) \bigr\}. 
\]
For any set $A\subset\Z^d$ denote by $A^c=\Z^d\backslash A$ its
complement and
by $\tau(A)$ the hitting time of $A$ by the random walk $(X_s)$,
and we write $\tau(z)$ for $\tau(\{z\})$ for any point $z\in\Z^d$.
Let us decompose the solution $u$ into $u=u_1+u_2$ according to the two
groups of paths (I) and (II)
mentioned in the \hyperref[intro]{Introduction}
\begin{eqnarray*}
u_1(t,z)&=&\E_0 \biggl[\exp \biggl\{\int
_0^t\xi(X_s)\,\dd s \biggr\}\id\{
X_t=z\} \id \bigl\{\tau\bigl(Z_t^{\ssup1}\bigr)\le
t,\tau \bigl(B_t^c \bigr)>t \bigr\} \biggr],
\\
u_2(t,z)&=&\E_0 \biggl[\exp \biggl\{\int
_0^t\xi(X_s)\,\dd s \biggr\}\id\{
X_t=z\} \id \bigl\{\tau\bigl(Z_t^{\ssup1}\bigr)> t
\mbox{ or } \tau \bigl(B_t^c \bigr)\le t \bigr\} \biggr].
\end{eqnarray*}

In Lemma~\ref{ll_u2} below, we use the results from Section~\ref{s_neg}
to prove that the total mass of $u_2$ is negligible.
In order to prove that $u_1$ localises around $Z_t^{\ssup1}$, we
introduce the gap
\[
\mathfrak{g}_t=\xi\bigl(Z_t^{\ssup1}\bigr)-\max
\bigl\{\xi(z)\dvtx z\in B_t \backslash\bigl\{Z_t^{\ssup1}
\bigr\} \bigr\}
\]
between the value of the potential $\xi$ at the point $Z_t^{\ssup1}$
and in the
rest of the ball $B_t$. In Lemma~\ref{l_gap}
we find a lower bound for $\mathfrak{g}_t$. This bound tends to
infinity for $\gamma<1$ but is going to
zero for $1\le\gamma<2$. However, the lower bound turns out to be just
large enough
to provide localisation of the principal eigenfunction of the
Anderson Hamiltonian $\Delta+\xi$ around $Z_t^{\ssup1}$, which is
proved in
Lemma~\ref{l_eig}.
This easily implies the localisation of $u_1$ around $Z_t^{\ssup1}$ and
allows us to prove Theorem~\ref{main_w} in the end of this section.

%
\begin{lemma}
\label{ll_u2}
Almost surely,
\[
\biggl\{U(t)^{-1}\sum_{z\in\Z^d}u_2(t,z)
\biggr\} \1_{\mathcal{E}_c(t)}\to0 \qquad\mbox{as }t\to\infty.
\]
\end{lemma}

\begin{pf}
We have
%
%
\begin{equation}
\label{ind}\quad  \sum_{z\in\Z^d}u_2(t,z) =
\E_0 \biggl[\exp \biggl\{\int_0^t
\xi(X_s)\,\dd s \biggr\}\id \bigl\{\tau\bigl(Z_t^{\ssup1}
\bigr)> t\mbox{ or } \tau \bigl(B_t^c \bigr)\le t \bigr\}
\biggr].
\end{equation}
Observe that if a path belongs to the set in the indicator function
above then either it passes through $Z_t^{\ssup1}$ and reaches
the maximum of the potential there but leaves the ball $B_t$ thus
belonging to $E_1(t)$, or it reaches the maximum of the
potential not in $Z_t^{\ssup1}$ thus belonging to $E_2(t)$ or $E_3(t)$,
depending on whether the maximum of the potential over
the path exceeds the value $\xi_{r_tg_t}^{\ssup{k_t}}$. Hence, we have
on the event $\mathcal{E}_c(t)$
\begin{eqnarray*}
\sum_{z\in\Z^d}u_2(t,z) &\le&
\E_0 \biggl[\exp \biggl\{\int_0^t
\xi(X_s)\,\dd s \biggr\}\id_{E_1(t)\cup
E_2(t)\cup E_3(t)} \biggr]
\\
&=&U_1(t)+U_2(t)+U_3(t).
\end{eqnarray*}
The statement of the lemma now follows from Lemmas~\ref{l_u1}, \ref
{l_u2} and~\ref{l_u3}.
\end{pf}

%
\begin{lemma}
\label{l_gap}
On the event $\mathcal{E}_c(t)$, the gap $\mathfrak{g_t}$ is positive and,
for any $\varepsilon>0$,
\[
\log\mathfrak{g}_t> (1/\gamma-1-\varepsilon)\log\log t
\]
eventually for all $t$.
\end{lemma}

\begin{pf} Let $z\in B_t\backslash\{Z_t^{\ssup1}\}$.
Then $\varPsi_t(z)\le\varPsi_t(Z_t^{\ssup2})$ and we have on the
event $\mathcal{E}_c(t)$
\begin{eqnarray*}
d_t\lambda_t&\le&\varPsi_t
\bigl(Z_t^{\ssup1}\bigr)- \varPsi_t
\bigl(Z_t^{\ssup2}\bigr)  \le \varPsi_t
\bigl(Z_t^{\ssup1}\bigr)-\varPsi_t(z)
\\
&=&\xi\bigl(Z_t^{\ssup1}\bigr)-\xi(z)+\frac{|z|-|Z_t^{\ssup1}|}{\gamma
t}\log
\log t. 
\end{eqnarray*}
Since $|Z_t^{\ssup1}|<r_tg_t$ on the event $\mathcal{E}_c(t)$, the
last term satisfies
\[
\frac{|z|-|Z_t^{\ssup1}|}{\gamma t}\log\log t\le\frac{|Z_t^{\ssup
1}|\rho_t}{\gamma
t}\log\log t\le
\frac{r_tg_t\rho_t}{\gamma t}\log\log t=O(d_tg_t\rho_t).
\]
We obtain uniformly for all $z\in B_t\backslash\{Z_t^{\ssup1}\}$
\[
d_t\lambda_t\le\xi\bigl(Z_t^{\ssup1}
\bigr)-\xi(z)+O(d_tg_t \rho_t)
\]
and so
\[
\mathfrak{g}_t\ge d_t\lambda_t+O(d_tg_t
\rho_t)=d_t\lambda_t+o(d_t
\lambda_t)
\]
on according to \eqref{b}. This estimate implies the statement of the lemma
since $\log d_t\sim(\frac{1}{\gamma}-1)\log\log t$ and $\lambda_t$ is
negligible according to \eqref{a}.
\end{pf}

Let $\gamma_t$ and $v_t$ be the principal eigenvalue and eigenfunction
of $\Delta+\xi$
with zero boundary conditions in the ball $B_t$. We extend $v_t$ by
zero to the whole space $\Z^d$ and
we assume that $v_t$ is normalised so that $v_t(Z_t^{\ssup1})=1$. The
eigenfunction $v_t$ has the following
probabilistic representation
\[
v_t(z)=\E_z \biggl[\exp \biggl\{\int
_0^{\tau(Z_t^{\ssup1})} \bigl(\xi(X_s)-
\gamma_t \bigr)\,ds \biggr\}\id \bigl\{\tau\bigl(Z_t^{\ssup1}
\bigr)<\tau \bigl( \Z^d\backslash B_t \bigr) \bigr\}
\biggr].
\]

%
\begin{lemma}
\label{l_eig}
Almost surely,
\[
\biggl\{\Vert v_t\Vert_2^2 \sum
_{z\in B_t\backslash\{Z_t^{\ssup1}\}} v_t(z) \biggr\}\1 _{\mathcal
{E}_c(t)}\to0
\qquad\mbox{as }t\to\infty.
\]
\end{lemma}

\begin{pf}
Consider the event $\mathcal{E}_c(t)$ and suppose that $t$ is
sufficiently large.
For each $n, p\in\N$ and $z\in B_t\backslash\{Z_t^{\ssup1}\}$ denote
\[
\mathcal{P}_{n,p}(t,z)= \bigl\{y\in\mathcal{P}_n\dvtx
y_0=z, y_n=Z_t^{\ssup1},
y_{i}\in B_t \backslash Z_t^{\ssup1}\
\forall i< n, p(y)=p \bigr\}.
\]

Integrating with respect to the waiting times $(\tau_i)$ of the random
walk, which are independent and
exponentially distributed with parameter $2d$
and observing that the probability of the first $n$ steps of the random
walk to follow a given geometric
path is $(2d)^{-n}$ we get
\begin{eqnarray*}
v_t(z)&=&\sum_{n\ge|z-Z_t^{\ssup1}|}\sum
_{p\le n}\sum_{y\in\mathcal
{P}_{n,p}(t,z)}(2d)^{-n}
{\mathsf{E}} \Biggl[\exp \Biggl\{\sum_{i=0}^{n-1}
\bigl(\xi(y_i)-\gamma_t \bigr)\tau_i \Biggr
\} \Biggr]
\\
&=&\sum_{n\ge|z-Z_t^{\ssup1}|}\sum_{1\le p\le n}
\sum_{y\in\mathcal{P}_{n,p}(t,z)} \prod_{i=0}^{n-1}
\int_0^{\infty} 
\exp \bigl\{- \bigl(
\gamma_t+2d-\xi(y_i) \bigr)t \bigr\}\,\dd t.
\end{eqnarray*}
The Rayleigh--Ritz formula implies
\begin{eqnarray*}
\gamma_t &=&\sup \bigl\{ \bigl\langle(\Delta+\xi)\phi, \phi \bigr
\rangle\dvtx\phi\in\ell^2(B_t), \phi|_{\partial B_t}=0,
\Vert\phi\Vert_2=1 \bigr\}
\\
&\ge& \bigl\langle(\Delta+\xi)\id_{\{Z_t^{\ssup1}\}}, \id_{\{
Z_t^{\ssup1}\}} \bigr
\rangle=\xi\bigl(Z_t^{\ssup1}\bigr)-2d
\end{eqnarray*}
and so for all $i$
%
%
\begin{equation}
\label{gg} \gamma_t+2d-\xi(y_i)\ge\xi
\bigl(Z_t^{\ssup1}\bigr)- \xi(y_i)\ge
\mathfrak{g}_t.
\end{equation}
Since $\mathfrak{g}_t>0$ eventually on the event $\mathcal{E}_c(t)$ by
Lemma~\ref{l_gap}, we use \eqref{gg}
to compute
%
%
\begin{eqnarray}
\label{ooo} v_t(z)&=&\sum_{n=|z-Z_t^{\ssup1}|}^{\infty}
\sum_{p\le n}\sum_{y\in\mathcal
{P}_{n,p}(t,z)}
\prod_{i=0}^{n-1}\frac{1}{\gamma_t+2d-\xi(y_i)}
\nonumber
\\
&\le&\sum_{p\ge|z-Z_t^{\ssup1}|}\sum_{n\ge p}
\sum_{y\in\mathcal{P}_{n,p}(t,z)} \prod_{i=0}^{n-1}
\frac{1}{\xi(Z_t^{\ssup1})-\xi(y_i)}
\nonumber
\\[-8pt]
\\[-8pt]
&\le&\sum_{p\ge|z-Z_t^{\ssup1}|}\sum_{n\ge p}(2d)^{-n}
\max_{y\in\mathcal
{P}_{n,p}(t,z)} \Biggl\{(2d)^{2n}\prod
_{i=0}^{n-1}\frac{1}{\xi(Z_t^{\ssup1})-\xi(y_i)} \Biggr\}
\nonumber
\\
&\le&\sum_{p\ge|z-Z_t^{\ssup1}|} \exp\max_{n\ge p}
\max_{y\in\mathcal
{P}_{n,p}(t,z)} \Biggl\{2n\log(2d)-\sum
_{i=0}^{n-1}\log \bigl(\xi\bigl(Z_t^{\ssup1}
\bigr)-\xi(y_i) \bigr) \Biggr\}
\nonumber
\end{eqnarray}
since $\sum_{n\ge p}(2d)^{-n}\le1$ for $p\ge1$.
%
Fix\vspace*{1pt} some positive $\varepsilon\in(\frac{1} \gamma-1, \frac{1}{\gamma}-\frac{1} 2
)$. Notice that this is possible since $\gamma<2$ and so
$\frac{1} \gamma-\frac{1} 2>0$.
Let $p\ge|z-Z_t^{\ssup1}|$, $n\ge p$, and $y\in\mathcal{P}_{n,p}(t,z)$.
By Lemma~\ref{app}(e), the set $G_{r_tg_t}$ is totally disconnected
and so
%
%
\begin{equation}
\label{ooo1} n_{+}(y,G_{r_tg_t})\le \biggl\lceil
\frac{n+1}{2} \biggr\rceil\le\frac{n}{2}+1. 
\end{equation}
In each point $y_i\in G_{r_tg_t}$, we can estimate by Lemma~\ref{l_gap}
%
%
\begin{equation}
\label{ooo2} \log \bigl(\xi\bigl(Z_t^{\ssup1}\bigr)-
\xi(y_i) \bigr)\ge\log \mathfrak{g}_t> (1/\gamma-1-
\varepsilon)\log\log t.
\end{equation}
On the other hand,
%
%
\begin{equation}
\label{ooo3} n_{-}(y,G_{r_tg_t})=n+1-n_{+}(y,G_{r_tg_t})
\ge\frac{n}{2}
\end{equation}
and in each point $y_i\notin G_{r_tg_t}$ we get by Lemma~\ref{app}(d)
%
%
\begin{equation}
\label{ooo4} \log \bigl(\xi\bigl(Z_t^{\ssup1}\bigr)-
\xi(y_i) \bigr)=\log \bigl( \xi_{r_tg_t}^{\ssup
{k_t}}-
\xi_{r_tg_t}^{\ssup{m_t}} \bigr) >(1/ \gamma-\varepsilon)\log\log t
\end{equation}
%
by Lemma~\ref{l_gap}.
Using \eqref{ooo2} and \eqref{ooo4} and taking into account
that the last point $Z_t^{\ssup1}$ of the path belongs to $G_{r_tg_t}$
but does not contribute to the sum, we obtain
\begin{eqnarray*}
&&2n\log(2d)-\sum_{i=0}^{n-1}\log \bigl(
\xi\bigl(Z_t^{\ssup1}\bigr)-\xi(y_i) \bigr)
\\
&&\qquad\le2n\log(2d)- \bigl(n_{+}(y,G_{r_tg_t})-1 \bigr) (1/
\gamma-1-\varepsilon)\log\log t
\\
&&\qquad\phantom{\le}-n_{-}(y,G_{r_tg_t}) (1/\gamma-
\varepsilon)\log\log t.
\end{eqnarray*}
Since $\frac{1} \gamma-1-\varepsilon<0$ and $\frac{1} \gamma
-\varepsilon>0$, we
can estimate further using \eqref{ooo1} and~\eqref{ooo3}
\begin{eqnarray*}
&&2n\log(2d)-\sum_{i=0}^{n-1}\log \bigl(
\xi\bigl(Z_t^{\ssup1}\bigr)-\xi(y_i) \bigr)
\\
&&\qquad\le2n\log(2d)-\frac{n}{2}(1/\gamma-1-\varepsilon)\log \log t -
\frac{n}{2}(1/\gamma-\varepsilon)\log\log t
\\
&&\qquad=2n\log(2d)-n(1/\gamma-1/2-\varepsilon)\log\log t.
\end{eqnarray*}
Since $\frac{1} \gamma-\frac{1} 2-\varepsilon>0$, this function is
decreasing in
$n$ and can be estimated
by its value at $n=p$. This implies
\begin{eqnarray*}
&&2n\log(2d)-\sum_{i=0}^{n-1}\log \bigl(
\xi\bigl(Z_t^{\ssup1}\bigr)-\xi(y_i) \bigr)
\\
&&\qquad\le2p\log(2d)-p(1/\gamma-1/2-\varepsilon)\log\log t \le -p\delta\log\log
t
\end{eqnarray*}
with some $\delta>0$. Substituting this into \eqref{ooo}, we obtain
\[
v_t(z)\le\sum_{p\ge|z-Z_t^{\ssup1}|}(\log
t)^{-p\delta} 
\le2(\log t)^{-\delta|z-Z_t^{\ssup1}|}.
\]
Since $v_t(z)$ decays geometrically in distance of $z$ from $Z_t^{\ssup1}$,
$(\log t)^{-\delta}\to0$, and $v_t(Z_t^{\ssup1})=1$, the statement
of the lemma is now obvious.
\end{pf}
%

%
\begin{remark}\label{rem2}
Observe that, similarly to the proof of Proposition~\ref{Hh}, we have a
competition of the positive and negative
terms in the sum in \eqref{ooo}, and we want the negative terms to
dominate. The contribution of the positive terms is
of order $(1/\gamma-1)\log\log t$ and the contribution of the negative
terms is roughly $(1/\gamma)\log\log t$.
This leads to the condition $1-1/\gamma<1/\gamma$, which restricts our
proof to the case $0<\gamma<2$.
\end{remark}

\begin{pf*}{Proof of Theorem~\ref{main_w}}
We have
%
%
\begin{eqnarray}
\label{aaa} 1-\frac{u(t,Z_t^{\ssup1})}{U(t)} &=&U(t)^{-1}\sum
_{z\neq Z_t^{\ssup1}}u(t,z)
\nonumber
\\[-8pt]
\\[-8pt]
&\le& U(t)^{-1}\sum_{z\neq Z_t^{\ssup1}}u_1(t,z)+U(t)^{-1}
\sum_{z\in\Z^d}u_2(t,z).
\nonumber
\end{eqnarray}
The second term converges to zero on the event $\mathcal{E}_c(t)$ by
Lemma~\ref{ll_u2}. The first term satisfies the conditions of \cite
{GKM06}, Theorem~4.1, with $B=B_t$, $V=\xi$, and $\varGamma=\{
Z_t^{\ssup1}\}$,
which implies that, for all $z\in B_t$,
\[
u_1(t,z)\le u_1\bigl(t,Z_t^{\ssup1}
\bigr)\Vert v_t\Vert _2^2 v_t(z).
\]
Observing that $U(t)\ge u_1(t,Z_t^{\ssup1})$ and $u_1(t,z)=0$ for
$z\notin B_t$,
we obtain
\[
U(t)^{-1}\sum_{z\neq Z_t^{\ssup1}}u_1(t,z)
\le\Vert v_t\Vert_2^2\sum
_{z\in
B_t\backslash
\{Z_t^{\ssup1}\}}v_t(z),
\]
which converges to zero on the event $\mathcal{E}_c(t)$ by Lemma~\ref
{l_eig}. As
both terms in \eqref{aaa} converge to zero on the event $\mathcal
{E}_c(t)$ and
$\Prob\{\mathcal{E}_c(t)\}\to1$ by
Proposition~\ref{l_4den}(b), we obtain that
\[
1-\frac{u(t,Z_t^{\ssup1})}{U(t)}\to0 \qquad\mbox{as }t\to\infty
\]
in probability.
\end{pf*}
%


\section{Ageing}
\label{s_age}

In this section, we discuss the ageing behaviour of the parabolic
Anderson model.
Throughout this section, we assume that $\gamma>0$. As we pointed out
in the \hyperref[intro]{Introduction},
although the results proved in this section hold for all $\gamma>0$,
they only imply ageing of the
parabolic Anderson model for $0<\gamma<2$ as otherwise the solution $u$
may not be localised at $Z_t^{\ssup1}$.

We begin by showing
that whenever the maximiser of $\varPsi$ has moved from one point to
another, it cannot go back to the original point.

%
\begin{lemma}
\label{a1}
For $s>0$, $\{ T_t>s\}= \{Z_t^{\ssup1}=Z_{t+s}^{\ssup1} \}$ eventually
for all $t$.
\end{lemma}

\begin{pf} If $T_t>s$, then $Z_t^{\ssup1}=Z_{t+s}^{\ssup1}$ by the
definition of $T_t$. Suppose $Z_t^{\ssup1}=Z_{t+s}^{\ssup1}$
but there is $u\in(t,t+s)$ such that $Z_t^{\ssup1}\neq Z_u^{\ssup
1}$. Consider
an auxiliary function $\phi\dvtx[t,t+s]\to\R$ given by
\[
\phi(x)=\varPsi_x\bigl(Z_t^{\ssup1}\bigr)-
\varPsi_x \bigl(Z_u^{\ssup1} \bigr)=\xi
\bigl(Z_t^{\ssup1}\bigr)-\xi \bigl(Z_u^{\ssup
1}
\bigr)-\frac{|Z_t^{\ssup1}|-|Z_u^{\ssup1}|}{\gamma x} \log\log x. 
\]
Observe that
\[
\phi'(x)=\frac{|Z_t^{\ssup1}|-|Z_u^{\ssup1}|}{\gamma x^2\log
x}(\log x\log\log x-1)
\]
and so $\phi'$ does not change the sign on the interval $[t,t+s]$ if
$t$ is large enough. Hence, $\phi$
is strictly monotone on $[t,t+s]$. However, this contradicts the
observation that
$\phi(t)\ge0$ (since $Z_t^{\ssup1}$ is the maximiser of $\varPsi_t$ and
$Z_u^{\ssup1}\neq Z_t^{\ssup1}$), $\phi(u)\le0$
(since $Z_u^{\ssup1}$ is the maximiser of $\varPsi_u$ and $Z_t^{\ssup
1}\neq
Z_u^{\ssup1}$), and $\phi(t+s)\ge0$
(since $Z_t^{\ssup1}=Z_{t+s}^{\ssup1}$ is the maximiser of $\varPsi
_{t+s}$ and
$Z_{u}^{\ssup1}\neq Z_t^{\ssup1}$).
\end{pf}

Now we are going to compute the probability of $ \{Z_t^{\ssup1}
=Z_{t+wt}^{\ssup1} \}$, $w>0$, using
the point processes $\varPi_t\equiv\varPi_{t,0}$ studied in
Section~\ref{s_ppp}.
However, we need to restrict them to a finite box growing to infinity
to justify integration and
passing to the limit. In order to do so, for
each $n\in\N$, we define the event
\[
\mathcal{A}(n,w,t)= \bigl\{Y_{t,Z_t^{\ssup1}}\ge-n, \varPsi _{t+wt}(z)
\le \varPsi_{t+wt}\bigl(Z_t^{\ssup1}\bigr)\ \forall z\in
\Z^d \mbox{ s.t. }Y_{t,z}\ge-n \bigr\}
\]
and show that $\Prob\{Z_t^{\ssup1}=Z_{t+wt}^{\ssup1} \}$ is captured by
the probabilities of these events.

%
\begin{lemma}
\label{a2}
For any $w>0$,
\[
\lim_{t\to\infty}\Prob \bigl\{Z_t^{\ssup1}=Z_{t+wt}^{\ssup1}
\bigr\}=\lim_{n\to
\infty}\lim_{t\to\infty}\Prob \bigl\{
\mathcal{A}(n,w,t) \bigr\},
\]
provided the limit on the right-hand side exists.
\end{lemma}

\begin{pf}
To obtain an upper bound, observe that
%
%
\begin{equation}
\label{b1} \Prob \bigl\{Z_t^{\ssup1}=Z_{t+wt}^{\ssup1}
\bigr\} \le \Prob \bigl\{\mathcal{A}(n,w,t) \bigr\} +\Prob\{Y_{t,Z_t^{\ssup
1}}
\le-n \}.
\end{equation}
By Proposition~\ref{l_4den},
%
%
\begin{equation}
\label{con1} \lim_{n\to\infty} \lim_{t\to\infty} \Prob
\{Y_{t,Z_t^{\ssup1}}\le-n \} =\lim_{n\to\infty}\Prob \bigl
\{Y^{\ssup1}\le-n \bigr\}=0.
\end{equation}

For a lower bound, we have
%
%
\begin{equation}
\label{b2} \Prob \bigl\{Z_t^{\ssup1}=Z_{t+wt}^{\ssup1}
\bigr\} \ge \Prob \bigl\{\mathcal{A}(n,w,t) \bigr\} -\Prob\{ Y_{t,Z_{t+wt}^{\ssup1}}
\le-n \}.
\end{equation}
Observe that for all $z$ we have, as $t\to\infty$,
%
%
\begin{eqnarray}
\label{psiss} \varPsi_{t+wt}(z)&=&\xi(z)-\frac{|z|}{\gamma
(t+wt)}\log\log
(t+wt)
\nonumber
\\
&=&\varPsi_t(z)+\frac{w|z|}{(1+w)\gamma t} \bigl(\log\log t +o(1) \bigr)
\\
&=&\varPsi_t(z)+d_{r_t}\frac{w\theta}{1+w}
\frac{|z|}{r_t} \bigl(1 +o(1) \bigr)
\nonumber
\end{eqnarray}
and so the condition $Y_{t,Z_{t+wt}^{\ssup1}}\le-n$ is equivalent to
%
%
\begin{equation}
\label{co1} \frac{\varPsi_{t+wt}(Z_{t+wt}^{\ssup
1})-a_{r_t}}{d_{r_t}}-\frac
{w\theta
}{1+w}\frac{|Z_{t+wt}^{\ssup1}|}{r_t}
\bigl(1+o(1) \bigr)\le-n.
\end{equation}
It is easy to see that $r_{t+wt}\sim(1+w)r_t$. This implies that
$d_{r_{t+wt}}\sim d_{r_t}$ and
\[
a_{r_{t+wt}}-a_{r_t}\sim d_{r_t}\gamma^{-1}d
\log(1+w).
\]
Now condition \eqref{co1} is equivalent to
\[
\biggl[\frac{\varPsi_{t+wt}(Z_{t+wt}^{\ssup
1})-a_{r_{t+wt}}}{d_{r_{t+wt}}} +\gamma^{-1}d\log(1+w)-w\theta
\frac{|Z_{t+wt}^{\ssup
1}|}{r_{t+wt}} \biggr] \bigl(1+o(1) \bigr)\le-n
\]
and by Proposition~\ref{l_4den} we obtain
%
%
\begin{eqnarray}
\label{con2} &&\lim_{n\to\infty} \lim_{t\to\infty}
\Prob\{Y_{t,Z_{t+wt}^{\ssup1}}\le-n \}
\nonumber
\\
&&\qquad=\lim_{n\to\infty}\lim_{t\to\infty}\Prob \biggl
\{ \biggl[Y_{t+wt,Z_{t+wt}^{\ssup1}}+\gamma^{-1}d\log(1+w) -w\theta
\frac{|Z_{t+wt}^{\ssup1}|}{r_{t+wt}} \biggr]
\nonumber
\\[-8pt]
\\[-8pt]
&&\phantom{\hspace*{230pt}} {}\times \bigl(1+o(1) \bigr) \le-n \biggr\}
\nonumber
\\
&&\qquad=\lim_{n\to\infty}\Prob
\bigl\{Y^{\ssup1}+\gamma^{-1}d\log(1+w)-w\theta\bigl |X^{\ssup1}\bigr |
\le-n \bigr\}=0.
\nonumber
\end{eqnarray}
Combining the bounds \eqref{b1} and \eqref{b2} with the convergence
results \eqref{con1} and \eqref{con2},
we obtain the required statement.
\end{pf}

Now we show that the probabilities of the events $\mathcal{A}(n,w,t)$
converge to a finite explicit integral.

%
\begin{lemma}
\label{a3}
For any $w\ge0$,
\[
\lim_{n\to\infty}\lim_{t\to\infty}\Prob \bigl\{\mathcal
{A}(n,w,t) \bigr\} =\int_{\R^d\times\R}\exp \bigl\{-\nu
\bigl(D_{w}(x,y) \bigr) \bigr\}\nu(dx,dy)<\infty,
\]
where $D_w(x,y)$ has been defined in \eqref{dc}.
\end{lemma}

\begin{pf}
We have
\begin{eqnarray*}
&&\Prob \bigl\{\mathcal{A}(n,w,t) \bigr\}
\\
&&\qquad=\int_{\R^d\times[-n,\infty)} \Prob\bigl\{ \bigl(Z_t^{\ssup
1}r_t^{-1},Y_{t,Z_t^{\ssup1}}
\bigr)\in\dd x\times\dd y,
\\
&&\phantom{\qquad=\int_{\R^d\times[-n,\infty)} \Prob \bigl\{}
\varPsi_{t+wt}(z)\le\varPsi_{t+wt}\bigl(Z_t^{\ssup1}
\bigr)\ \forall z\in \Z^d \mbox{ s.t. }Y_{t,z}\ge-n \bigr\}.
\end{eqnarray*}
Observe that according to \eqref{psiss}
the condition $\varPsi_{t+wt}(z)\le\varPsi_{t+wt}(Z_t^{\ssup1})$ is
equivalent to
\[
\varPsi_t(z)+d_{r_t}\frac{w\theta}{1+w}\frac{|z|}{r_t}
\bigl(1 +o(1) \bigr) \le\varPsi_t\bigl(Z_t^{\ssup1}
\bigr)+d_{r_t} \frac{w\theta}{1+w}\frac{|Z_t^{\ssup1}|}{r_t} \bigl(1 +o(1)
\bigr),
\]
that is, to
\[
Y_{t,z}+\frac{w\theta}{1+w}\frac{|z|}{r_t} \bigl(1 +o(1) \bigr) \le
Y_{t,Z_t^{\ssup1}}+\frac{w\theta}{1+w}\frac{|Z_t^{\ssup
1}|}{r_t} \bigl(1 +o(1) \bigr).
\]

Consider the point process $\varPi_t$ on $\hat H_{-n}^{-\alpha}$, where
$\alpha\in(\theta\frac{w}{1+w}, \theta)$.
The requirement
\[
\bigl\{ \bigl(Z_t^{\ssup1}r_t^{-1},Y_{t,Z_t^{\ssup1}}
\bigr)\in \dd x\times\dd y, \varPsi_{t+wt}(z)\le\varPsi_{t+wt}
\bigl(Z_t^{\ssup1}\bigr)\ \forall z\in\Z^d \mbox{
s.t. }Y_{t,z}\ge-n \bigr\}
\]
means that $\varPi_t$ has one point in $\dd x\times\dd y$ and no points
in the domain
\begin{eqnarray*}
D_{n, w,t}(x,y)&=& \bigl(\R^d\times[y,\infty ) \bigr)
\\
&&{}\cup \biggl\{(\bar x,\bar y)\in\R^d\times[-n,\infty)\dvtx
\\
&&\phantom{{}\cup \biggl\{}y+
\frac{w\theta|x|}{1+w} \bigl(1 +o(1) \bigr)\le\bar y+\frac{w\theta
|\bar
x|}{1+w} \bigl(1 +o(1)
\bigr) \biggr\}.
\end{eqnarray*}
Hence, by Lemma~\ref{l_ppp},
\begin{eqnarray*}
&&\lim_{t\to\infty}\Prob \bigl\{\mathcal{A}(n,w,t) \bigr\}
\\
&&\qquad=\lim_{t\to\infty}\int_{\R^d\times[-n,\infty)}\Prob \bigl
\{\varPi_t(\dd x\times\dd y)=1, \varPi_t
\bigl(D_{n, w,t}(x,y) \bigr)=0 \bigr\}
\\
&&\qquad= \int_{\R^d\times[-n,\infty)}\Prob \bigl\{\varPi(\dd x\times\dd y)=1,
\varPi \bigl(D_{n, w}(x,y) \bigr)=0 \bigr\}
\\
&&\qquad=\int_{\R^d\times[-n,\infty)}\exp \bigl\{-\nu \bigl(D_{n, w}(x,y)
\bigr) \bigr\}\nu(dx,dy),
\end{eqnarray*}
where
\[
D_{n, w}(x,y)=D_w(x,y)\cap \bigl(\R^d
\times[-n,\infty ) \bigr).
\]
Taking the limit in this way is justified as $\hat H_{-n}^{-\alpha}$ is
compact and contains $\R^d\times[-n,\infty)$.

It remains to show that
%
%
\begin{eqnarray}
\label{last} &&\lim_{n\to\infty}\int_{\R^d\times[-n,\infty
)}\exp
\bigl\{-\nu \bigl(D_{n,
w}(x,y) \bigr) \bigr\}\nu(dx,dy)
\nonumber
\\[-8pt]
\\[-8pt]
&&\qquad= \int_{\R^d\times\R}\exp \bigl\{-\nu \bigl(D_{w}(x,y)
\bigr) \bigr\}\nu(dx,dy)<\infty.
\nonumber
\end{eqnarray}

Observe that $\nu(D_{n,w}(x,y))\ge\nu(\R^d\times(y,\infty))$ for all
$x\in\R^d$ and $y\ge-n$. Then
\[
\id_{\R^d\times[-n,\infty)}(x,y)\exp \bigl\{-\nu \bigl(D_{n,
w}(x,y) \bigr)
\bigr\} \le\exp \bigl\{ -\nu \bigl(\R^d\times(y,\infty) \bigr) \bigr
\}.
\]
It is easy to see that $\exp\{-\nu(\R^d\times(y,\infty))\}$ is
integrable with respect to the measure $\nu$ on $\R^d\times\R$ since
using \eqref{dom} and the substitution $u=e^{-\gamma y}$ we get
%
%
\begin{eqnarray}
\label{domint} &&\int_{\R^d\times\R} \exp \bigl\{-\nu \bigl(
\R^d\times(y,\infty) \bigr) \bigr\}\nu(\dd x,\dd y)
\nonumber
\\
&&\qquad=\int_{-\infty}^{\infty}\int_{\R^d}
\gamma\exp \bigl\{-\gamma y -\gamma\theta|x|-2^d(\gamma
\theta)^{-d}e^{-\gamma y} \bigr\}\,\dd x \,\dd y
\nonumber
\\[-8pt]
\\[-8pt]
&&\qquad=2^d(\gamma\theta)^{-d}\int_{-\infty}^{\infty}
\gamma\exp \bigl\{ -\gamma y -2^d(\gamma\theta)^{-d}e^{-\gamma y}
\bigr\}\,\dd y
\nonumber
\\
&&\qquad=2^d(\gamma\theta)^{-d}\int_{0}^{\infty}
\exp \bigl\{ -2^d(\gamma\theta)^{-d}u \bigr\}\,\dd u=1.
\nonumber
\end{eqnarray}
Now \eqref{last} follows from the dominated convergence theorem.
\end{pf}

Finally, we combine all results of this section to prove ageing.

\begin{pf*}{Proof of Theorem~\ref{main_a}}
For any $w>0$, we have by Lemmas~\ref{a1}, \ref{a2} and~\ref{a3},
\begin{eqnarray*}
F(w)&:=&\lim_{t\to\infty}\Prob\{T_t/t\le w \} =1-\lim
_{t\to\infty}\Prob \bigl\{Z_t^{\ssup1}=Z_{t+wt}^{\ssup1}
\bigr\}
\\
&=&1-\lim_{n\to\infty}\lim_{t\to\infty}\Prob \bigl\{
\mathcal{A}(n,w,t) \bigr\}
\\
&=&1-\int_{\R^d\times\R}\exp \bigl\{-\nu \bigl(D_{w}(x,y)
\bigr) \bigr\}\nu(dx,dy).
\end{eqnarray*}
Observe that $\exp\{-\nu(D_{w}(x,y))\}\le\exp\{-\nu(\R^d\times
(y,\infty))\}$ which is integrable with respect to the measure
$\nu$ by \eqref{domint}. Since $\nu(D_{w}(x,y))\to\break \nu(D_{w_0}(x,y))$
whenever $w\to w_0\in(0,\infty)$ the function $F$
is continuous.

If $w\to0{\scriptscriptstyle+}$ then $\nu(D_{w}(x,y)\to\nu(\R
^d\times(y,\infty))$ and by \eqref{domint} we obtain
\[
\lim_{w\to0{\scriptscriptstyle+}}F(w)=1-\int_{\R^d\times\R}\exp \bigl\{-
\nu \bigl(\R^d\times(y,\infty) \bigr) \bigr\}\nu(dx,dy)=0.
\]

Finally, if $w\to\infty$ then $\nu(D_{w}(x,y))\to\nu(D_{\infty
}(x,y))$, where
\[
D_{\infty}(x,y)= \bigl\{(\bar x,\bar y)\in\R^d\times\R\dvtx y+
\theta|x|\le\bar y+\theta|\bar x| \bigr\} \cup \bigl(\R^d\times[y,
\infty ) \bigr).
\]
Compute
\begin{eqnarray*}
\nu \bigl(D_{\infty}(x,y) \bigr) 
&\ge&
\int_{|\bar x|> |x|}\int_{y+\theta|x|-\theta|\bar x|}^{\infty
}\gamma
\exp \bigl \{-\gamma\bar y -\gamma\theta|\bar x| \bigr \}\,\dd\bar y\,\dd\bar x
\\
&=&\exp\bigl \{-\gamma y-\gamma\theta|x| \bigr \}\int
_{|\bar x|>
|x|}\dd\bar x=\infty.
\end{eqnarray*}
Hence, $F(w)\to1$ as $w\to\infty$.
\end{pf*}



%

\printaddresses


\begin{thebibliography}{18}

\bibitem{An58}
%
\begin{barticle}[auto:STB|2014/02/12|12:18:25]
\bauthor{\bsnm{Anderson},~\bfnm{P.~W.}\binits{P.~W.}}
(\byear{1958}).
\btitle{Absence of diffusion in certain random lattices}.
\bjournal{Phys. Rev.}
\bvolume{109}
\bpages{1492--1505}.
\end{barticle}
%
\bptok{imsref}%
\endbibitem

\bibitem{AD}
%
\begin{barticle}[mr]
\bauthor{\bsnm{Aurzada},~\bfnm{Frank}\binits{F.}} \AND
\bauthor{\bsnm{D{\"o}ring},~\bfnm{Leif}\binits{L.}}
(\byear{2011}).
\btitle{Intermittency and ageing for the symbiotic branching model}.
\bjournal{Ann. Inst. Henri Poincar\'e Probab. Stat.}
\bvolume{47}
\bpages{376--394}.
\bid{doi={10.1214/09-AIHP356}, issn={0246-0203}, mr={2814415}}
\end{barticle}
%
\bptok{imsref}%
\endbibitem

\bibitem{BC}
%
\begin{bincollection}[mr]
\bauthor{\bsnm{Ben Arous},~\bfnm{G{\'e}rard}\binits{G.}} \AND
\bauthor{\bsnm{{\v{C}}ern{\'y}},~\bfnm{Ji{\v{r}}{\'{\i}}}\binits{J.}}
(\byear{2006}).
\btitle{Dynamics of trap models}.
In \bbooktitle{Mathematical Statistical Physics}
\bpages{331--394}.
\bpublisher{Elsevier},
\blocation{Amsterdam}.
\bid{doi={10.1016/S0924-8099(06)80045-4}, mr={2581889}}
\end{bincollection}
%
\bptok{imsref}%
\endbibitem

\bibitem{BK}
%
\begin{barticle}[mr]
\bauthor{\bsnm{Biskup},~\bfnm{Marek}\binits{M.}} \AND
\bauthor{\bsnm{K{\"o}nig},~\bfnm{Wolfgang}\binits{W.}}
(\byear{2001}).
\btitle{Long-time tails in the parabolic {A}nderson model with bounded
potential}.
\bjournal{Ann. Probab.}
\bvolume{29}
\bpages{636--682}.
\bid{doi={10.1214/aop/1008956688}, issn={0091-1798}, mr={1849173}}
\end{barticle}
%
\bptok{imsref}%
\endbibitem

\bibitem{CM94}
%
\begin{barticle}[mr]
\bauthor{\bsnm{Carmona},~\bfnm{Ren{\'e}~A.}\binits{R.~A.}} \AND
\bauthor{\bsnm{Molchanov},~\bfnm{S.~A.}\binits{S.~A.}}
(\byear{1994}).
\btitle{Parabolic {A}nderson problem and intermittency}.
\bjournal{Mem. Amer. Math. Soc.}
\bvolume{108}
\bpages{viii$+$125}.
\bid{doi={10.1090/memo/0518}, issn={0065-9266}, mr={1185878}}
\end{barticle}
%
\bptok{imsref}%
\endbibitem

\bibitem{DD}
%
\begin{barticle}[mr]
\bauthor{\bsnm{Dembo},~\bfnm{Amir}\binits{A.}} \AND
\bauthor{\bsnm{Deuschel},~\bfnm{Jean-Dominique}\binits{J.-D.}}
(\byear{2007}).
\btitle{Aging for interacting diffusion processes}.
\bjournal{Ann. Inst. Henri Poincar\'e Probab. Stat.}
\bvolume{43}
\bpages{461--480}.
\bid{doi={10.1016/j.anihpb.2006.07.001}, issn={0246-0203}, mr={2329512}}
\end{barticle}
%
\bptok{imsref}%
\endbibitem

\bibitem{GK05}
%
\begin{bincollection}[mr]
\bauthor{\bsnm{G{\"a}rtner},~\bfnm{J{\"u}rgen}\binits{J.}} \AND
\bauthor{\bsnm{K{\"o}nig},~\bfnm{Wolfgang}\binits{W.}}
(\byear{2005}).
\btitle{The parabolic {A}nderson model}.
In \bbooktitle{Interacting Stochastic Systems}
\bpages{153--179}.
\bpublisher{Springer},
\blocation{Berlin}.
\bid{doi={10.1007/3-540-27110-4_8}, mr={2118574}}
\end{bincollection}
%
\bptok{imsref}%
\endbibitem

\bibitem{GKM06}
%
\begin{barticle}[mr]
\bauthor{\bsnm{G{\"a}rtner},~\bfnm{J{\"u}rgen}\binits{J.}},
\bauthor{\bsnm{K{\"o}nig},~\bfnm{Wolfgang}\binits{W.}} \AND
\bauthor{\bsnm{Molchanov},~\bfnm{Stanislav}\binits{S.}}
(\byear{2007}).
\btitle{Geometric characterization of intermittency in the parabolic
{A}nderson model}.
\bjournal{Ann. Probab.}
\bvolume{35}
\bpages{439--499}.
\bid{doi={10.1214/009117906000000764}, issn={0091-1798}, mr={2308585}}
\end{barticle}
%
\bptok{imsref}%
\endbibitem

\bibitem{GM90}
%
\begin{barticle}[mr]
\bauthor{\bsnm{G{\"a}rtner},~\bfnm{J.}\binits{J.}} \AND
\bauthor{\bsnm{Molchanov},~\bfnm{S.~A.}\binits{S.~A.}}
(\byear{1990}).
\btitle{Parabolic problems for the {A}nderson model. {I}.
{I}ntermittency and related topics}.
\bjournal{Comm. Math. Phys.}
\bvolume{132}
\bpages{613--655}.
\bid{issn={0010-3616}, mr={1069840}}
\end{barticle}
%
\bptok{imsref}%
\endbibitem

\bibitem{GM98}
%
\begin{barticle}[mr]
\bauthor{\bsnm{G{\"a}rtner},~\bfnm{J.}\binits{J.}} \AND
\bauthor{\bsnm{Molchanov},~\bfnm{S.~A.}\binits{S.~A.}}
(\byear{1998}).
\btitle{Parabolic problems for the {A}nderson model. {II}.
{S}econd-order asymptotics and structure of high peaks}.
\bjournal{Probab. Theory Related Fields}
\bvolume{111}
\bpages{17--55}.
\bid{doi={10.1007/s004400050161}, issn={0178-8051}, mr={1626766}}
\end{barticle}
%
\bptok{imsref}%
\endbibitem

\bibitem{GS}
%
\begin{barticle}[mr]
\bauthor{\bsnm{G{\"a}rtner},~\bfnm{J{\"u}rgen}\binits{J.}} \AND
\bauthor{\bsnm{Schnitzler},~\bfnm{Adrian}\binits{A.}}
(\byear{2011}).
\btitle{Time correlations for the parabolic {A}nderson model}.
\bjournal{Electron. J. Probab.}
\bvolume{16}
\bpages{1519--1548}.
\bid{doi={10.1214/EJP.v16-917}, issn={1083-6489}, mr={2827469}}
\end{barticle}
%
\bptok{imsref}%
\endbibitem

\bibitem{KLMS09}
%
\begin{barticle}[mr]
\bauthor{\bsnm{K{\"o}nig},~\bfnm{Wolfgang}\binits{W.}},
\bauthor{\bsnm{Lacoin},~\bfnm{Hubert}\binits{H.}},
\bauthor{\bsnm{M{\"o}rters},~\bfnm{Peter}\binits{P.}} \AND
\bauthor{\bsnm{Sidorova},~\bfnm{Nadia}\binits{N.}}
(\byear{2009}).
\btitle{A two cities theorem for the parabolic {A}nderson model}.
\bjournal{Ann. Probab.}
\bvolume{37}
\bpages{347--392}.
\bid{doi={10.1214/08-AOP405}, issn={0091-1798}, mr={2489168}}
\end{barticle}
%
\bptok{imsref}%
\endbibitem

\bibitem{LM}
%
\begin{bmisc}[auto:STB|2014/02/12|12:18:25]
\bauthor{\bsnm{Lacoin},~\bfnm{H.}\binits{H.}} \AND
\bauthor{\bsnm{M{\"o}rters},~\bfnm{P.}\binits{P.}}
(\byear{2012}).
\bhowpublished{A scaling limit theorem for the parabolic Anderson model
with exponential potential.
In \textit{Probability in Complex Physical Systems.
In Honour of J. G\"artner and E. Bolthausen}.
\textit{Springer Proc.  Math.} \textbf{11} 247--271. Springer, Berlin.}
\end{bmisc}
%
\bptok{imsref}%
\endbibitem

\bibitem{M94}
%
\begin{bincollection}[mr]
\bauthor{\bsnm{Molchanov},~\bfnm{S.}\binits{S.}}
(\byear{1994}).
\btitle{Lectures on random media}.
In \bbooktitle{Lectures on Probability Theory ({S}aint-{F}lour, 1992)}.
\bseries{Lecture Notes in Math.}
\bvolume{1581}
\bpages{242--411}.
\bpublisher{Springer},
\blocation{Berlin}.
\bid{doi={10.1007/BFb0073874}, mr={1307415}}
\end{bincollection}
%
\bptok{imsref}%
\endbibitem

\bibitem{MOS}
%
\begin{barticle}[mr]
\bauthor{\bsnm{M{\"o}rters},~\bfnm{Peter}\binits{P.}},
\bauthor{\bsnm{Ortgiese},~\bfnm{Marcel}\binits{M.}} \AND
\bauthor{\bsnm{Sidorova},~\bfnm{Nadia}\binits{N.}}
(\byear{2011}).
\btitle{Ageing in the parabolic {A}nderson model}.
\bjournal{Ann. Inst. Henri Poincar\'e Probab. Stat.}
\bvolume{47}
\bpages{969--1000}.
\bid{doi={10.1214/10-AIHP394}, issn={0246-0203}, mr={2884220}}
\end{barticle}
%
\bptok{imsref}%
\endbibitem

\bibitem{HKM06}
%
\begin{barticle}[mr]
\bauthor{\bsnm{van~der Hofstad},~\bfnm{Remco}\binits{R.}},
\bauthor{\bsnm{K{\"o}nig},~\bfnm{Wolfgang}\binits{W.}} \AND
\bauthor{\bsnm{M{\"o}rters},~\bfnm{Peter}\binits{P.}}
(\byear{2006}).
\btitle{The universality classes in the parabolic {A}nderson model}.
\bjournal{Comm. Math. Phys.}
\bvolume{267}
\bpages{307--353}.
\bid{doi={10.1007/s00220-006-0075-4}, issn={0010-3616}, mr={2249772}}
\end{barticle}
%
\bptok{imsref}%
\endbibitem

\bibitem{HMS08}
%
\begin{barticle}[mr]
\bauthor{\bsnm{van~der Hofstad},~\bfnm{Remco}\binits{R.}},
\bauthor{\bsnm{M{\"o}rters},~\bfnm{Peter}\binits{P.}} \AND
\bauthor{\bsnm{Sidorova},~\bfnm{Nadia}\binits{N.}}
(\byear{2008}).
\btitle{Weak and almost sure limits for the parabolic {A}nderson model
with heavy tailed potentials}.
\bjournal{Ann. Appl. Probab.}
\bvolume{18}
\bpages{2450--2494}.
\bid{doi={10.1214/08-AAP526}, issn={1050-5164}, mr={2474543}}
\end{barticle}
%
\bptok{imsref}%
\endbibitem

\bibitem{ZM87}
%
\begin{barticle}[mr]
\bauthor{\bsnm{Zel'dovich},~\bfnm{Ya.~B.}\binits{Ya.~B.}},
\bauthor{\bsnm{Molchanov},~\bfnm{S.~A.}\binits{S.~A.}},
\bauthor{\bsnm{Ruzma{\u\i}kin},~\bfnm{A.~A.}\binits{A.~A.}} \AND
\bauthor{\bsnm{Sokolov},~\bfnm{D.~D.}\binits{D.~D.}}
(\byear{1987}).
\btitle{Intermittency in random media}.
\bjournal{Uspekhi Fiz. Nauk}
\bvolume{152}
\bpages{3--32}.
\bid{doi={10.1070/PU1987v030n05ABEH002867}, mr={0921018}}
\end{barticle}
%
\bptok{imsref}%
\endbibitem

\end{thebibliography}
\end{document}